\documentclass[letterpaper]{article} 
\usepackage{aaai2026}  
\usepackage{amsfonts,amsmath}
\usepackage{times}  
\usepackage{helvet}  
\usepackage{courier}  
\usepackage{ntheorem}
\usepackage{cases}
\usepackage{array}
\usepackage[hyphens]{url}  
\usepackage{graphicx} 
\usepackage{natbib}  
\usepackage{caption} 
\usepackage{algorithm}
\usepackage{algorithmic}

\usepackage{multirow,booktabs,makecell}
\usepackage{newfloat}
\usepackage{listings}
\DeclareCaptionStyle{ruled}{labelfont=normalfont,labelsep=colon,strut=off} 
\floatstyle{ruled}
\newfloat{listing}{tb}{lst}{}
\floatname{listing}{Listing}
\newtheorem*{proof}{proof}

\newtheorem{Theorem}{Theorem}
\newtheorem{df}{Definition}
\newtheorem{Lemma}{Lemma}

\newtheorem{Remark}{Remark}

\newtheorem{Ass}{Assumption}

\title{S-D-RSM: Stochastic Distributed Regularized Splitting Method for Large-Scale Convex Optimization Problems}

\author {
    Maoran Wang\textsuperscript{\rm 1},
 Xingju Cai\textsuperscript{\rm 1,\rm 2},
 Yongxin Chen\textsuperscript{\rm 3}\thanks{Corresponding author}
}
\affiliations {
    \textsuperscript{\rm 1}School of Mathematical Sciences, Nanjing Normal University, P.R. China\\
    \textsuperscript{\rm 2}Key Laboratory of NSLSCS (NNU), Ministry of Education, P.R. China \\
    \textsuperscript{\rm 3}School of Mathematics and Statistics, Nanjing University of Science and Technology, P.R. China  \\
    \{230901010, caixingju\}@njnu.edu.cn, chenyongxin@buaa.edu.cn
}

\usepackage{bibentry}

\begin{document}

\maketitle

\begin{abstract}
This paper investigates  problems of large-scale distributed composite  convex optimization, with motivations from a broad range of applications, including multi-agent systems, federated learning, smart grids, wireless sensor networks, compressed sensing, and so on. Stochastic gradient descent (SGD) and its variants are commonly employed to solve such problems. However, existing algorithms often rely on vanishing step sizes, strong convexity assumptions, or entail substantial computational overhead to ensure convergence or obtain favorable complexity. To bridge the gap between theory and practice, we integrate consensus optimization and operator splitting techniques (see Problem Reformulation) to develop a novel stochastic splitting algorithm, termed the stochastic distributed regularized splitting method (S-D-RSM). In practice, S-D-RSM performs parallel updates of proximal mappings and gradient information for only a randomly selected subset of agents at each iteration. By introducing regularization terms, it effectively mitigates consensus discrepancies among distributed nodes. In contrast to conventional stochastic methods, our theoretical analysis establishes that S-D-RSM achieves global convergence without requiring diminishing step sizes or strong convexity assumptions. Furthermore, it achieves an iteration complexity of 1/epsilon with respect to both the objective function value and the consensus error. Numerical experiments show that S-D-RSM achieves up to two to three times speedup compared with state-of-the-art baselines, while maintaining comparable or better accuracy. These results not only validate the algorithm's theoretical guarantees but also demonstrate its effectiveness in practical tasks such as compressed sensing and empirical risk minimization. 
\end{abstract}

\section{Introduction}
In this work, we study the problem  large-scale distributed composite convex optimization:
\begin{equation}\label{p1}
  \min_{x\in \mathbb{R}^n} \left\{\Phi(x):=\sum_{i=1}^{m}\left(f_i(x)+g_i(x)\right)\right\},
\end{equation}
where $m$ is the numbers of nodes, $\{f_i\}_{i=1}^m$ is a sequence of proper, lower semicontinuous convex functions (not necessarily differentiable), $\{g_i\}_{i=1}^m$ is a sequence of convex functions that are Fréchet differentiable on $\mathbb{R}^n$, and each gradient $\nabla g_i$ is $\frac{1}{\beta_i}$-Lipschitz continuous. Throughout this paper, the usual restrictive requirement of strongly convexity of $f_i$ or $g_i$ is not needed \cite{pathak2020fedsplit,li2020secure,NEURIPS2024_e0b6f389,sadiev2024stochasticproximalpointmethods}. Problem~\eqref{p1} arises in a wide range of applications, including economics and traffic theory~\cite{2006Optimization,GU2019190}, image processing~\cite{2016An,ehrhardt2025guide}, machine learning~\cite{JMLR:v25:23-1040}, and other fields. 

Some ``full participation'' optimization methods—where all nodes are involved in computation at per iteration—have been proposed to solve problem \eqref{p1}; see, for example,~\cite{2011Generalized,2015Forward,pathak2020fedsplit,2021Distributed}. Although these algorithms admit global convergence under general convexity assumptions, their per-iteration cost remains high due to the need to compute all proximal mappings ${\rm prox}_{f_i}$  and evaluate all gradients $\nabla g_i$ for large-scale problems. As a result, stochastic (i.e., partial participation) optimization methods have attracted increasing attention.
\subsection{Related Works}
\paragraph{Gradient-Based Methods for Smooth Problems.}  
{Stochastic gradient descent (SGD)}~\cite{1951A} is a foundational algorithm widely used in machine learning. Since its introduction by Robbins and Monro \shortcite{1951A}, SGD has undergone numerous developments, giving rise to variants such as stochastic batch gradient descent~\cite{doi:10.1137/070704277} and compressed gradient descent~\cite{NIPS2017_6c340f25}. Recently, Gower et al.~\shortcite{pmlr-v97-qian19b} proposed a general framework for analyzing SGD with arbitrary sampling strategies in the strongly convex setting. Overall, while these methods are effective in many practical settings, their theoretical convergence guarantees typically rely on restrictive conditions such as vanishing step sizes and strong convexity.
\paragraph{Proximal Point Algorithms for Non-Smooth Problems.} 
For non-smooth optimization problems, the proximal point algorithm (PPA)~\cite{doi:10.1137/0314056} and its variants have been extensively investigated. Compared to  gradient-based methods, PPA exhibits greater robustness to inaccuracies in step size selection, as evidenced by the analyses in~\cite{ryu2014stochastic,parikh2014proximal}. For large-scale non-smooth  problems, stochastic variants of PPA (S-PPA) are more commonly employed in practice \cite{bertsekas2011incremental,sppa,patrascu2018nonasymptotic}. Under random sampling of component functions, vanishing step sizes, and suitable measurability and boundedness assumptions, Bianchi~\shortcite{sppa} established the almost sure convergence of S-PPA in the ergodic sense. Recently, Li et al.~\shortcite{NEURIPS2024_e0b6f389} proposed an extrapolated version of S-PPA (also known as FedExProx) for federated learning, which incorporates mini-batch sampling and an extrapolation step to accelerate convergence. Under convexity, Lipschitz continuity, and interpolation regimes~\cite{montanari2022interpolation}—which are satisfied in overparameterized deep learning models—they established an iteration complexity of $\mathcal{O}(\epsilon^{-1})$. Sadiev et al.~\shortcite{sadiev2024stochasticproximalpointmethods} further established linear convergence under additional structural assumptions.
\paragraph{Proximal Gradient Methods for Composite Problems.} 
For composite problems with \textit{multiple} smooth components and a \textit{single} non-smooth convex function, the stochastic proximal gradient (S-PG) method—originating from a combination of SGD~\cite{1951A} and proximal gradient methods~\cite{Beck2009}—has been extensively investigated~\cite{rosasco2016stochastic,atchade2017perturbed,spg}. In~\cite{rosasco2016stochastic}, the almost sure convergence of S-PG was established under strong convexity and vanishing step sizes. Under the general convex setting, Atchad\'{e} et al.~\shortcite{atchade2017perturbed} developed a unified analytical framework for both unbiased and biased gradient estimators in S-PG and derived an $\mathcal{O}(\epsilon^{-2})$ complexity bound under the assumption that the non-smooth component of the objective function is nonnegative and vanishing step sizes. Recently, Rosasco et al.~\shortcite{spg} refined the complexity result of Atchad\'{e} et al.~\shortcite{atchade2017perturbed} by establishing an improved bound of $\mathcal{O}\left(\epsilon^{1/(t - 1)}\right)$ without requiring the non-negativity assumption, with a vanishing step size of the form $\mathcal{O}(1/k^t)$, where $t \in (1/2,1)$ and $k$ denotes the iteration index. However, no convergence guarantees for S-PG are available in the absence of either strong convexity or vanishing step sizes.
\paragraph{Operator Splitting Methods.} When both smooth and non-smooth components are present in \textit{multiple} blocks, \textit{operator splitting techniques} provide a powerful algorithmic framework for designing deterministic algorithms. Some methods have been successfully extended to stochastic settings in recent works. Cevher et al.~\shortcite{NIPS2016_5d6646aa,Cevher2018} introduced the {stochastic forward Douglas-Rachford (S-FDR)} splitting method, establishing a stochastic extension of deterministic FDR~\cite{2015Forward}. Although S-FDR adopts the SGD-style gradient estimate, its requirement to compute all proximal mappings per-iteration raises scalability concerns for large-scale problems. Furthermore, by inheriting SGD's framework, S-FDR inherits similar theoretical requirements including vanishing step sizes and strong convexity assumptions~\cite{NIPS2016_5d6646aa}. More recently, a broader algorithmic framework was introduced by Combettes et al.~\cite{2015STOCHASTIC,9747479,2025gf}, who developed the stochastic generalized forward–backward (S-GFB) method as a stochastic extension of the deterministic GFB \cite{2011Generalized}. In contrast to S-FDR, this approach reduces the per-iteration computational burden by updating only a subset of the proximal mappings \( \{ \mathrm{prox}_{\gamma f_i} \}_{i=1}^m \), though it may still pose computational challenges due to the need for full gradient evaluations, particularly in large-scale applications.

For clarity and ease of comparison, the properties of the aforementioned algorithms are summarized in Table~\ref{t1}. In the table, $S_k$ denotes the index set sampled at iteration~$k$. The notions of ``Convergence'' and ``Complexity'' are analyzed under standard convexity assumptions, without imposing strong convexity or diminishing step-size conditions.
\begin{table*}[h]
\centering
\small
{\begin{tabular}{c|c|c|c|c|c}
\hline
{Algorithm} & $f_i,g_i$ & Step size  & Computational cost& Convergence & Complexity\\
\hline
\makecell{S-PPA\\ \cite{sppa}}& $g_i=0$  &vanishing & \({\rm prox}_{\gamma_k f_{i_k}}\) & No & No \\
\hline
\makecell{FedExProx\\ \cite{NEURIPS2024_e0b6f389}} & $g_i=0$ & constant & \(\left\{{\rm prox}_{\gamma f_i}\right\}_{i\in S_k}\)& No & No \\
\hline
\makecell{S-PG\\ \cite{spg}} & $f_i=0,~i\ge 2$ & vanishing  &  \(\{\nabla g_i\}_{i\in S_k}\) & No & \makecell{$\mathcal{O}(\epsilon^{1/(t- 1)})$\\$t\in(1/2,1)$}\\
\hline
\makecell{S-GFB\\ \cite{2025gf}} & $f_i,g_i \neq 0$ &  constant  & \(\{{\rm prox}_{\gamma f_{i}}\}_{i\in S_k}\); \(\{\nabla g_i\}_{i\in [m]}\) & Yes & No\\
\hline
\makecell{S-FDR\\  \cite{Cevher2018}}& $f_i,g_i \neq 0$ & vanishing  &  \(\{{\rm prox}_{\gamma f_{i}}\}_{i\in [m]}\); \(\{\nabla g_i\}_{i\in S_k}\) & Yes & No\\
\hline
\pmb{This paper} & \pmb{$f_i,g_i \neq 0$} &  \pmb{constant}& \(\{{\rm prox}_{\gamma f_{i}},\nabla g_i\}_{i\in {\pmb {S_k}}}\)  & \pmb{Yes} & \pmb{$\mathcal{O}{(\epsilon^{-1})}$}\\
\hline
\end{tabular}}
\caption{Comparison of the properties of S-D-RSM (Algorithm~\ref{A}) and several state-of-the-art methods. }
\label{t1}
\end{table*}
\paragraph{Theoretical and Practical Trade-offs.}
Compared to full participation approaches, stochastic methods significantly reduce computational costs by involving only a subset of nodes in each iteration. However, they also entail inherent trade-offs in step size policies, strong convexity requirements, gradient approximation accuracy, convergence guarantees, and complexity analysis. Vanishing step sizes are commonly used to establish almost sure convergence~\cite{sppa,spg}, but they may degrade practical performance. Conversely, constant step sizes typically require more accurate gradient estimates~\cite{2025gf,9747479,2015STOCHASTIC}, which can increase computational overhead. Furthermore, complexity analysis often relies on additional structural assumptions about the objective function, such as strong convexity or interpolation regimes~\cite{NIPS2016_5d6646aa,pmlr-v97-qian19b,spg,NEURIPS2024_e0b6f389}. These challenges highlight the need for algorithms that are both theoretically robust and computationally efficient.

\paragraph{Contributions.}Motivated by the unresolved theoretical-practical trade-offs in stochastic optimization, we propose a novel framework addressing three persistent limitations in state-of-the-art methods:

\begin{itemize}
    \item \pmb{Practical Limitations of Vanishing Step Sizes} \\
    While vanishing step sizes ($\gamma_k \to 0$) ensure theoretical convergence, empirical evidence consistently highlights adverse effects: asymptotic slowdown preventing $\epsilon$-optimal solutions and hyperparameter sensitivity causing sharp convergence deterioration \cite{doi:10.1137/16M1080173}. This bottleneck acutely impacts cross-device federated learning with heterogeneous compute capabilities.
    
    \item \pmb{Restrictive Functional Assumptions} \\
    Existing $\mathcal{O}(\epsilon^{-1})$ guarantees rely on structurally convenient but impractical conditions: strong convexity violated by large  machine learning models, and interpolation regimes implausible under non-IID data (e.g., recommendation systems) \cite{zhang2024goal}. 
    
    \item \pmb{Large-scale computing bottlenecks} \\
    Despite their stochastic formulations, prevalent operator-splitting methods still inherit deterministic burdens: S-FDR-type algorithms require $\mathcal{O}(m)$ proximal evaluations per iteration, while S-GFB-type methods necessitate $\mathcal{O}(m)$ gradient computations. These computational demands can become prohibitive in large-scale distributed systems, especially when $m$ is large or when proximal or gradient evaluations are costly.
\end{itemize}

The main contributions are summarized as follows:
\begin{itemize}
\item[$\bullet$] The proposed method integrates \textit{consensus optimization} with \textit{operator splitting} and exploits \textit{parallelism} by evaluating only a \textit{subset} of proximal mappings \( \{ \mathrm{prox}_{\gamma f_i} \}_{i=1}^m \) and gradients \( \{ \nabla g_i \}_{i=1}^m \) at each iteration. In addition, regularization is introduced into each subproblem to mitigate consensus discrepancies among distributed nodes, as confirmed by numerical experiments.
\item[$\bullet$] We provide a rigorous convergence analysis showing that S-D-RSM achieves global convergence under general convexity assumptions, \textit{without} requiring strong convexity, interpolation, or vanishing step sizes. The method attains a sublinear ergodic convergence rate of $\mathcal{O}(1/K)$ with respect to both the objective gap and consensus violation, leading to an  iteration complexity of $\mathcal{O}(\epsilon^{-1})$. Notably, we establish almost sure convergence of the iterate sequence, further reinforcing the algorithm’s reliability in practice.
\item[$\bullet$] Since the theoretical guarantees of our algorithm are established solely based on the objective function and do not depend on the underlying data distribution across devices, it retains global convergence and an $\mathcal{O}(\epsilon^{-1})$ complexity under heterogeneous settings, provided that the loss function is convex.
\end{itemize}
\section{Notations}\label{s2}
We denote by $\Gamma_0(\mathbb{R}^n)$ the set of all proper, lower semicontinuous, and convex functions on $\mathbb{R}^n$. 
Given $f \in \Gamma_0(\mathbb{R}^n)$, the subdifferential of $f$ is defined as
\[
\partial f : x \mapsto \left\{ u \in \mathbb{R}^n \mid f(y) \ge f(x) + \langle u, y - x \rangle,\ \forall y \in \mathbb{R}^n \right\},
\]
and its proximal mapping is defined by
\[
{\rm prox}_f : x \mapsto \arg\min_{y \in \mathbb{R}^n} \left\{ f(y) + \frac{1}{2} \| y - x \|^2 \right\}.
\]
From the definition of ${\rm prox}_f$, it can be verified that for all $u,x \in \mathbb{R}^n$ and $\delta > 0$,
\begin{equation}\label{xs}
x = {\rm prox}_{f}(u - \delta x) \Leftrightarrow x = {\rm prox}_{\frac{ f}{1+\delta}} \left( \frac{u}{1+\delta} \right).
\end{equation}
Let $(\Omega, \mathcal{F}, \mathbb{P})$ denote a probability space, with $\Omega$ the sample space, $\mathcal{F}$ the $\sigma$-algebra, and $\mathbb{P}$ the probability measure. The abbreviation ``$\mathbb{P}$-a.s.'' refers to ``$\mathbb{P}$-almost surely''. A sequence of random variables $\{\xi^k\}_{k=1}^\infty$ is said to converge $\mathbb{P}$-a.s. to a random variable $\xi$, if
\[
\mathbb{P}\left(\left\{\omega \in \Omega \mid \lim\limits_{k \to \infty} \xi^k(\omega) = \xi(\omega)\right\}\right) = 1,
\]
which is denoted as $\lim_{k\to\infty} \xi^k = \xi$, $\mathbb{P}$-a.s. Unless otherwise specified, all inequalities involving random variables are understood to hold $\mathbb{P}$-almost surely. The bold symbol $\boldsymbol{x}$ represents a vector of $m-1$ stacked vectors, $\{x_i\}_{i=1}^{m-1} \subset \mathbb{R}^n$, i.e., $\boldsymbol{x} = (x_1, x_2, \dots, x_{m-1})$. Moreover, we define $\sigma(x^0, x^1, \dots, x^k) \subset \mathcal{F}$ as the smallest $\sigma$-algebra generated by the set of random variables $\{x^0, x^1, \dots, x^k\}$. For a  random variable \( v \) and a \(\sigma\)-algebra \( \mathcal{J} \subset \mathcal{F} \), we denote by \( \mathbb{E}(v | \mathcal{J}) \) the conditional expectation of \( v \) given \( \mathcal{J} \), and write \( v \perp\perp \mathcal{J} \) to denote that \( v \) is independent of \( \mathcal{J} \). Finally, for any real number \( r \), the largest integer not greater than \( r \) is denoted by \( \lfloor r \rfloor \), and we define \( \frac{r}{0} = \infty \) in this paper.
\section{Problem Reformulations and the Proposed Algorithm}
In this section, we introduce two reformulations of problem \eqref{p1} that serve as the foundation of our approach. The first reformulation characterizes a system of equations satisfied by the solutions of problem~\eqref{p1}, offering theoretical guidance for the algorithm design and global convergence analysis. 
The second reformulation gives rise to the definition of $\epsilon$-optimal solutions, laying the foundation for the complexity analysis of the proposed algorithm. Since both reformulations are equivalent to problem~\eqref{p1}, their interrelationship is further clarified in Lemma~\ref{s12}.  Due to space constraints, all technical details and proofs are provided in the extended version.
\subsection{Problem Reformulation I}
\begin{Ass}
Assume that problem~\eqref{p1} admits at least one solution and problem~(\ref{p1}) satisfies
\[
 \bigcap_{i=1}^m \mathrm{ri}\left( \mathrm{dom} f_i \right)\neq\emptyset,
\]
where ``$\mathrm{ri}$'' denotes the set of relative interior points.
\end{Ass}

Under Assumption~1, we obtain \cite{1970ca}
\begin{equation}\label{cc}
{ \arg\min_{x \in \mathbb{R}^n}  \Phi(x) =\mathrm{zer}\left(\sum_{i=1}^{m}(\partial f_i+\nabla g_i)\right).}
\end{equation}
Based on~(\ref{cc}), we derive an alternative reformulation of the solution set of problem~(\ref{p1}).
\begin{Lemma}[Reformulation I]\label{kkt}
Let $\mathcal{S}$ be the set of all $(z_1, z_2, \ldots, z_{m-1}, x)$ satisfying the following system:
$$
\begin{cases}
\hspace{-0.8pt} x \hspace{-0.8pt} =\hspace{-0.8pt}  \operatorname{prox}_{\frac{\gamma f_m}{m-1}}\hspace{-4pt}\left(\hspace{-3.2pt}\tfrac{1}{m-1}\hspace{-2.2pt}\sum\limits_{i=1}^{m-1}\hspace{-3.4pt}\big(z_i \hspace{-1.1pt} - \hspace{-1.1pt}\sigma\gamma \nabla g_i(x)\big)\hspace{-1.1pt} -\hspace{-1.1pt} \frac{\gamma}{m-1} \nabla g_m(x)\hspace{-1.8pt} \right)\hspace{-.4pt},\\
\hspace{-0.8pt} x \hspace{-0.8pt} = \hspace{-0.8pt} \operatorname{prox}_{\gamma f_1}\hspace{-3.2pt}\left(2x - z_1 - (1-\sigma)\gamma \nabla g_1(x)\right), \\
\quad \vdots \\
\hspace{-0.8pt} x\hspace{-0.8pt}  = \hspace{-0.8pt} \operatorname{prox}_{\gamma f_{m-1}}\hspace{-3.4pt}\left(2x - z_{m-1} - (1-\sigma)\gamma \nabla g_{m-1}(x)\right).
\end{cases}
$$
Then the following assertions hold:
\begin{itemize}
\item[{\small$\bullet$}] If $x^\star$ minimizes $\Phi$, then there exist $z_1^\star, \ldots, z_{m-1}^\star \in \mathbb{R}^n$ such that $(z_1^\star, \ldots, z_{m-1}^\star, x^\star) \in \mathcal{S}$.
\item[{\small$\bullet$}]  Conversely, if $(z_1^\star, \ldots, z_{m-1}^\star, x^\star) \in \mathcal{S}$, then $x^\star$ minimizes $\Phi$.
\end{itemize}
\end{Lemma}
\subsection{Problem Reformulation II}
By introducing the constraint $ x_1 = x_2 =  \cdots = x_{m} $, problem \eqref{A} can be reformulated as:
\begin{align}
 &\min_{x_i}\sum_{i=1}^{m}\hspace{-1pt}f_i(x_i) + \hspace{-2pt}\sum_{i=1}^{m-1} \hspace{-2pt}\left\{(1{-}\sigma)g_i(x_i)\hspace{-1pt} +\hspace{-1pt} \sigma g_i(x_m)\right\} + g_m(x_m)\nonumber\\
 &~~{\rm s.t.}~ ~x_1 = x_2 = \cdots = x_{m}.\label{pd}
\end{align}
By leveraging the equivalence between problem \eqref{p1} and problem \eqref{pd}, the $\epsilon$-optimal solution of problem (\ref{p1}) is defined as follows.
\begin{df}\label{d1}
Let the tuple \( (x_1, x_2, \dots, x_m) \) consist of random variables generated by a stochastic algorithm over a probability space \( (\Omega, \mathcal{F}, \mathbb{P}) \). The tuple \( (x_1, x_2, \dots, x_m) \) is said to be \emph{\( \epsilon \)-optimal in expectation} if the random variables \( (x_1, x_2 \dots, x_{m}) \) satisfy the following two conditions for all \( i, j \in [m] \)
\[
\|\mathbb{E}[x_j - x_i]\| \leq \epsilon \quad \text{and} \quad \left|\mathbb{E}\left[H(x_1, \dots, x_{m})\right] - \Phi^\star\right| \leq \epsilon,
\]
where \( H \) denotes the objective function of the reformulated problem~\eqref{pd}, and \( \Phi^\star \) denotes the global optimal value of problem~\eqref{p1}.
\end{df}
\subsection{The Proposed Algorithm}
Based on the problem reformulation I, we introduce the stochastic distributed regularized splitting method (S-D-RSM) for addressing (1).
\begin{algorithm}[H]
\caption{S-D-RSM for solving problem (\ref{p1})}
\label{A}
\textbf{Input}: $K>0$; $\alpha_i\ge0$; $\sigma\in[0,1]$; $\alpha_i+\lfloor1-\sigma\rfloor\neq0$; initial point $ y_i^0, z_i^0\in\mathbb{R}^n$, $i\in[m-1]$. \\
\textbf{Parameter}:  \hspace{-2pt} $\gamma\hspace{-2pt}\in\hspace{-2pt}\left(0,\min\limits_{i\in[m-1]}\left\{\frac{2\alpha_i}{\frac{1}{\beta_m(m-1)}+\frac{\sigma}{\beta_i}},\frac{2(2+\alpha_i)\beta_i}{1-\sigma}\right\}\right)$; $\lambda_i\in\left(0,2+\alpha_i-\frac{(1-\sigma)\gamma}{2\beta_i}\right)$; error $=1$ and the error tolerance $\varepsilon>0$.\\
\textbf{Output}: Approximate solution $x^k$.\\
\textbf{Process}:
\begin{algorithmic}[1]
\STATE Let $k=0$, error $=1$.
\WHILE{error $>\varepsilon$ or $k\le K$}
    \STATE Server update
 { \begin{align*}
   x^{k+1}&\hspace{-3pt}= {\rm prox}_{\frac{\gamma f_m}{m-1}}\hspace{-1.5pt}\Bigg(\frac{1}{m-1}\hspace{-1pt}\sum_{i=1}^{m-1}\hspace{-1pt}\Big(z_i^k + \alpha_i(y_i^k-x^{k+1}) \\    &\hspace{-3pt}\quad -\frac{\gamma}{m-1}\hspace{-1pt}\nabla g_m(y_i^k)\Big) \hspace{-1pt}- \hspace{-1pt}\frac{\sigma\gamma}{m-1}\hspace{-1pt}\sum_{i=1}^{m-1}\hspace{-1.5pt}\nabla g_i(y_i^k)\Bigg).
    \end{align*}}
    
    \STATE Randomly select users $S_k \subseteq [m-1].$
    
    \STATE For user $i \in S_k$, compute
  {  \begin{align*}
 y_i^{k+1}&= {\rm prox}_{\gamma f_i}\Big(2x^{k+1}-z_i^{k} + \alpha_i(x^{k+1}-y_i^{k+1}) \\
    &\quad - (1-\sigma)\gamma\nabla g_i(x^{k+1})\Big), \\
    z_i^{k+1}&= z_i^{k} + \lambda_i\left(y_i^{k+1}-x^{k+1}\right).
    \end{align*}}
    
    \STATE For user $i \notin S_k$, set
  { \begin{align*}
  y_i^{k+1} &= y_i^{k}, \\
    z_i^{k+1} &= z_i^{k}.
    \end{align*}}
    \STATE Update error $\leftarrow \frac{\sum_{i=1}^{m-1}\|y_i^k-x^k\|^2}{\|x^k\|^2}$ and $k \leftarrow k+1$.
\ENDWHILE
\end{algorithmic}
\end{algorithm}
\begin{Remark}\label{remark1}
\begin{itemize}
\item[{\small$\bullet$}] Based on the definition of the proximal mapping, the subproblems in step 3 and step 5 contain regularization terms, specifically $\frac{1}{2\gamma}\sum_{i=1}^{m-1} \alpha_i\|x - y_i^k\|^2$ and $\frac{\alpha_i}{2\gamma}\|y_i - x^{k+1}\|^2$, which are introduced to balance the discrepancy between $x^{k+1}$ and $y_i^k$, with the balancing strength controlled by the parameter $\alpha_i$.
\item[{\small$\bullet$}] For computational convenience, we explicitly express $x^{k+1}$ and $y_i^{k+1}$ (for $i \in S_k$) based on~(\ref{xs}) as follows:
\begin{equation*}
{ {\begin{cases}
\begin{aligned}
\hspace{-2.0pt}x^{k+1}\hspace{-5.2pt}&=\hspace{-1.5pt}\operatorname{prox}_{\tfrac{\gamma f_m}{(1+\bar{\alpha})(m-1)}}
\hspace{-2.8pt}\Big(\tfrac{1}{(m-1)(1+\bar{\alpha})}\sum_{i=1}^{m-1}\Big(z_i^k + \alpha_i y_i^k \\
&\quad -\tfrac{\gamma}{m-1}\nabla g_m(y_i^k)\Big) \hspace{-1.4pt}-\hspace{-1.4pt} \tfrac{\sigma\gamma}{(m-1)(1+\bar{\alpha})}\hspace{-2.7pt}\sum_{i=1}^{m-1}\nabla g_i(y_i^k)\Big), \\
\hspace{-2.0pt}y_i^{k+1}\hspace{-5.2pt}&=\hspace{-1.5pt}\operatorname{prox}_{\tfrac{\gamma f_i}{1+\alpha_i}}\hspace{-5.4pt}\Big(\hspace{-2.2pt}\tfrac{2+\alpha_i}{1+\alpha_i}x^{k+1}\hspace{-8.5pt} -\hspace{-0.5pt} \tfrac{z_i^k}{1+\alpha_i} \hspace{-2pt}-\tfrac{(1-\sigma)\gamma}{1+\alpha_i}\hspace{-1.2pt}\nabla\hspace{-2pt} g_i(x^{k+1})\hspace{-2.2pt}\Big), 
\end{aligned}
\end{cases}}}
\end{equation*}
 where $\bar{\alpha}=\frac{1}{m-1}\sum_{i=1}^{m-1}\alpha_i.$ 
\item[$\bullet$] If \( \sigma > 0 \), the computation of \( x^{k+1} \) requires all \( \nabla g_i(y_i^k) \) and \( \nabla g_m(y_i^k) \) only at the initial iteration \( k = 0 \), while for \( k > 0 \), only a subset of these gradients needs to be computed. 
  
 \end{itemize}
\end{Remark} 
 \begin{Ass}
    Select $S_k\subset [m-1]$ such that $S_k\perp\perp\mathcal{F}_k$  with $\mathbb{P}(i\in S_k)=p_i>0, i\in[m-1]$, where the $\sigma$-algebra \( \mathcal{F}_k \) is defined as  
\[
\mathcal{F}_k = \sigma \left( \left\{ \boldsymbol{y}^j, \boldsymbol{z}^j \right\}_{j=0}^k \right),
\]  
where \( \boldsymbol{y}^k = \left( y_i^k \right)_{i=1}^{m-1} \) and \( \boldsymbol{z}^k = \left( z_i^k \right)_{i=1}^{m-1} \).
\end{Ass}
Consequently, the iterates \( x^{k} \), \( \boldsymbol{y}^k \), and \( \boldsymbol{z}^k \) are random variables, and the random set \( S_k \) is independent of the history \( \left\{ \boldsymbol{y}^j, \boldsymbol{z}^j \right\}_{j \leq k} \) and \( \left\{ x^j \right\}_{j \leq k+1} \). Furthermore, since \( x^{k+1} \) is generated through a continuous mapping of \( \boldsymbol{y}^k \) and \( \boldsymbol{z}^k \), it follows that \( x^{k+1} \) is \( \mathcal{F}_k \)-measurable.
\section{Main Theory Results}
In this section, we present several convergence results for Algorithm~\ref{A}. All theoretical results concerning Algorithm~\ref{A} are derived under Assumptions~1–2 and the parameter settings for $\gamma$, $\sigma$, $\alpha_i$, and $\lambda_i$ as specified in Algorithm~\ref{A}.  
To facilitate the convergence analysis of Algorithm~\ref{A}, we introduce the following auxiliary variables, which are not computed in practice:
 \begin{equation}\label{al3}
 \begin{cases}
 \begin{aligned}
\tilde{y}_i^{k+1}\hspace{-2pt}=& {\rm prox}_{\gamma f_i}\Big(2x^{k+1}-z_i^{k} + \alpha_i(x^{k+1}-\tilde{y}_i^{k+1}) \\
 & - (1-\sigma)\gamma\nabla g_i(x^{k+1})\Big),~ \forall i\in[m-1],\\
\tilde{z}_i^{k+1}\hspace{-2pt}=& z_i^{k} + \lambda_i\left(\tilde{y}_i^{k+1}-x^{k+1}\right),~ \forall i\in[m-1].
 \end{aligned}
 \end{cases}
 \end{equation}
 
 The following result demonstrates the decreasing properties of the random variables generated by Algorithm \ref{A}.
\begin{Lemma}[Decreasing properties]
{ The random sequence $\left\{ x^k, ( y_i^k,z_i^k)_{i=1}^{m-1} \right\}_{k=0}^\infty$ generated by Algorithm \ref{A} and the virtual user variables \( \left\{ (\tilde{y}_i^k)_{i=1}^{m-1} \right\}_{k=1}^\infty \) defined by \eqref{al3} satisfy that
\begin{equation}\label{e5}
\begin{aligned}
& \mathbb{E}\left\{ \sum_{i=1}^{m-1} \left( \frac{1}{\lambda_i p_i} \| z_i^{k+1} - z_i^\star \|^2 + \frac{\alpha_i}{p_i} \| y_i^{k+1} - x^\star \|^2\right) \bigg| \mathcal{F}_k \right\} \\
\le & \sum_{i=1}^{m-1}\left( \frac{1}{\lambda_i p_i} \| z_i^k - z_i^\star \|^2 +  \frac{\alpha_i}{p_i} \| y_i^k - x^\star \|^2 \right)\\
& - \sum_{i=1}^{m-1} \left( \alpha_i - \frac{\gamma}{2 \beta_m (m-1)} - \frac{\sigma \gamma}{2 \beta_i} \right) \| x^{k+1} - y_i^k \|^2 \\
& - \sum_{i=1}^{m-1} \left( 2 + \alpha_i - \frac{(1 - \sigma) \gamma}{2 \beta_i} - \lambda_i \right) \| x^{k+1} - \tilde{y}_i^{k+1} \|^2,
\end{aligned}
\end{equation}
for any \( (z_1^\star, z_2^\star, \dots, z_{m-1}^\star, x^\star) \in \mathcal{S} \).}
\end{Lemma}

Based on the reformulation I (Lemma 1) and the decreasing properties of Algorithm \ref{A} (Lemma 2), we are now able to demonstrate the \textit{\bf global convergence} of the sequence produced by Algorithm \ref{A}  {\bf without} requiring diminishing step sizes or strong convexity assumptions.
\begin{Theorem}[Convergence]\label{tc1}
Let $\left\{ x^k, ( y_i^k,z_i^k)_{i=1}^{m-1} \right\}_{k=0}^\infty$  denote the sequence generated by Algorithm \ref{A}. Then, the following hold:
\begin{itemize}
\item[{\small$\bullet$}] $\lim\limits_{k \to \infty} \| x^{k+1} - y_i^k \|=\lim\limits_{k \to \infty} \| x^k - \tilde{y}_i^k \|=0$, $\mathbb{P}$-a.s., $\forall i\in[m-1]$.
\item[{\small$\bullet$}] There exists a random variable $\tilde{x}$ taking values in $\arg\min_{x \in \mathbb{R}^n} \left\{ \Phi(x) \right\}$ such that $\lim_{k\to\infty}x^k =\tilde{x}$, $\mathbb{P}$-a.s. 
\end{itemize}
\end{Theorem} 

The following lemma establishes a connection between two equivalent formulations of problem~\eqref{p1}. Specifically, it presents properties of the objective function in problem~\eqref{pd} associated with the solution set $\mathcal{S}$ defined in Lemma~\ref{kkt}.
\begin{Lemma}\label{s12}
For any $(z_1^\star,z_2^\star,\dots,z_{m-1}^\star,x^\star) \in \mathcal{S}$ and any $x_i \in \mathbb{R}^n$ for $i \in [m]$, the following inequality holds:
\begin{equation}\label{efr1}
\small H(x_1,\cdots,x_m)-\Phi^\star\ge \frac{1}{\gamma}\sum_{i=1}^{m-1}\left\langle x^\star - z_i^\star, x_m - x_i \right\rangle,
\end{equation} 
where $\Phi^\star$ is the global optimal value of problem~\eqref{p1}.
\end{Lemma}

Next, we analyze the evolution of the objective function of problem \eqref{p1} along the sequence of iterates generated by Algorithm \ref{A}.
\begin{Lemma}
Let $\left\{ x^k, ( y_i^k,z_i^k)_{i=1}^{m-1} \right\}_{k=0}^\infty$ be the sequence generated by {Algorithm \ref{A}}, and let the virtual variables  \( \left\{ (\tilde{y}_i^k)_{i=1}^{m-1} \right\}_{k=1}^\infty \) be defined by \eqref{al3}. Then, for any $(z_1^\star,z_2^\star,\dots,z_{m-1}^\star,x^\star)\in \mathcal{S}$, the following inequality holds:
\begin{equation}\label{efr2}
\begin{aligned}
&2\gamma\left(H(\tilde{y}_1^{k+1},\cdots,\tilde{y}_{m-1}^{k+1}, x^{k+1})-H(x^\star,\cdots,x^\star)\right)\\
\leq&2\sum_{i=1}^{m-1}\left\langle x^\star-z_i^\star,x^{k+1}-\tilde{y}_i^{k+1}\right\rangle+ a_k-\mathbb{E}[a_{k+1}\mid\mathcal{F}_k]\\
&-\sum_{i=1}^{m-1}\left(2+\alpha_i-\lambda_i-\frac{(1-\sigma)\gamma}{ \beta_i}\right)\|\tilde{y}_i^{k+1}-x^{k+1}\|^2\\
&-\sum_{i=1}^{m-1}\left(\alpha_i-\frac{\gamma }{(m-1)\beta_m}-\frac{\sigma\gamma}{\beta_i}\right)\|x^{k+1}-y_i^{k}\|^2,
\end{aligned}
\end{equation}

where   
$$
a_k=\sum_{i=1}^{m-1}\left(\frac{\alpha_i}{p_i}\|y_i^k-x^\star\|^2+\frac{1}{\lambda p_i}\|z_i^k-z_i^\star\|^2\right).
$$
\end{Lemma}

Building on the previously established descent properties of Algorithm~\ref{A} and the structural characteristics of the original problem, we now establish the convergence rate of Algorithm~\ref{A} under a {\bf constant} step size and {\bf general convexity} assumptions.
\begin{Theorem}[Rate]\label{com}
Let \( \left\{ (z_i^k)_{i=1}^{m-1},  x^k, (y_i^k)_{i=1}^{m-1} \right\}_{k=0}^\infty \) be the sequence generated by Algorithm \ref{A}. Then for every $K \in \mathbb{N}$  and $i \in [m-1]$, define
\[
x_{\mathrm{av}}^K = \frac{1}{K} \sum_{k=0}^{K-1}x^{k+1},~ y_{\mathrm{av},i}^K = \frac{1}{K} \sum_{k=0}^{K-1}y_i^{k}.
\]
Then the following hold:
\begin{itemize}
\item[{\small$\bullet$}] $\left\| \mathbb{E}\left[x_{\mathrm{av}}^K - {y}_{\mathrm{av},i}^K\right]\right\| = \mathcal{O}\left({1}/{K}\right)$, for all $i \in [m-1]$.
\item[{\small$\bullet$}] $\left|\mathbb{E}\left[H({y}_{{\rm av},1}^K,\cdots,y_{ {\rm av},m-1}^K,x^K_{\rm av})\right]-\Phi^\star\right|=\mathcal{O}\left(1/{K}\right).$
\end{itemize}
\end{Theorem}

As a consequence of Theorem~\ref{com}, Algorithm~\ref{A} achieves an $\epsilon$-optimal solution in expectation within at most $\mathcal{O}(\epsilon^{-1})$ iterations.
\section{Numerical Experiments}
In this section, we apply the proposed S-D-RSM  to solve the compressed sensing problem and the logistic regression problem. We compare the performance of S-D-RSM with three state-of-the-art methods: Split-Douglas-Rachford method (SDR)\cite{2021Split}, S-GFB \cite{2025gf}, and S-FDR \cite{Cevher2018}. All algorithms were implemented in MATLAB 2021b, and experiments were conducted on a desktop computer equipped with an Intel(R) Core(TM) i5-10210U CPU @ 1.60GHz, 2112 MHz, and 8 GB RAM.

For all tested algorithms, each numerical experiment is repeated 20 times, and the average performance is reported. In each iteration, 30\% of the users are activated according to a uniform sampling strategy in the stochastic method. The initial points are set to zero vectors and are identical across all algorithms. The regularization parameters \(\alpha_i\), for all \(i \in [m-1]\), are set to 1, and the parameter \(\sigma\) is set to \(1/2\) in the proposed S-D-RSM algorithm. To ensure a fair comparison, the parameters of each algorithm are tuned as large as possible while still guaranteeing convergence.
\subsection{Compressed Sensing}
We begin by evaluating the empirical performance of the proposed method on the compressed sensing problem:
\begin{equation}\label{cs}
  \begin{aligned}
  &\min_{x\in\mathbb{R}^n}~\|x\|_1\\
  &~\text{s.t.}~Ax=b,
  \end{aligned}
\end{equation}
where $A \in \mathbb{R}^{p \times n}$ is the sensing matrix and $b \in \mathbb{R}^p$ is the observed measurement vector. Let $a_i$ denote the $i$-th row of $A$, and $b_i$ denote the $i$-th entry of $b$. By incorporating indicator functions for affine constraints, which equal 0 on the constraint set and $\infty$ otherwise, problem \eqref{cs} can be reformulated in the structure of problem~\eqref{p1}, with $g_i \equiv 0$:
\begin{equation}\label{cs1}
\min_{x\in\mathbb{R}^n}~\|x\|_1+\sum_{i=1}^{p}I_{C_i}(x),
\end{equation}
where each $C_i = \{x \mid a_i^\top x = b_i\}$ is a affine constraint set. 

\begin{figure}[htbp]
\begin{minipage}{0.49\columnwidth}
 \includegraphics[scale=0.42]{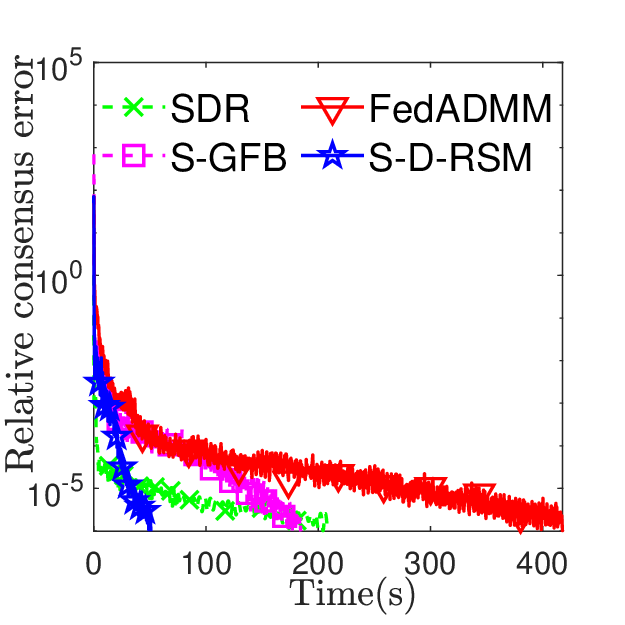} 
\end{minipage}
\begin{minipage}{0.49\columnwidth}
 \includegraphics[scale=0.42]{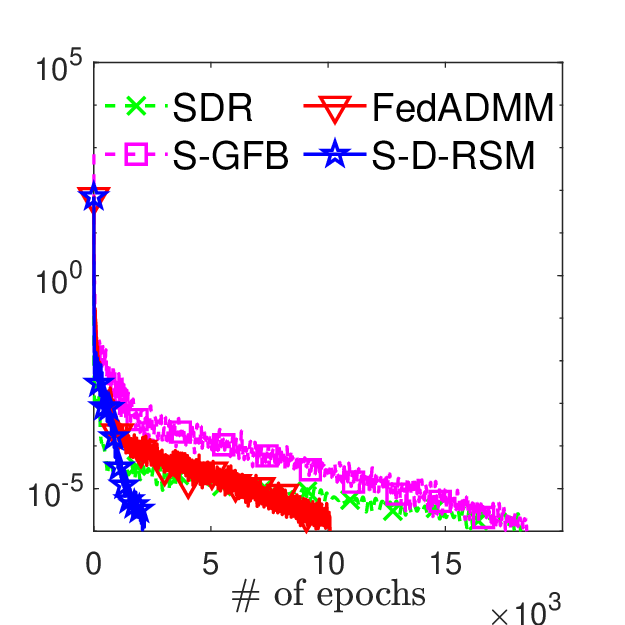} 
 \end{minipage}\\
 \begin{minipage}{0.49\columnwidth}
 \includegraphics[scale=0.42]{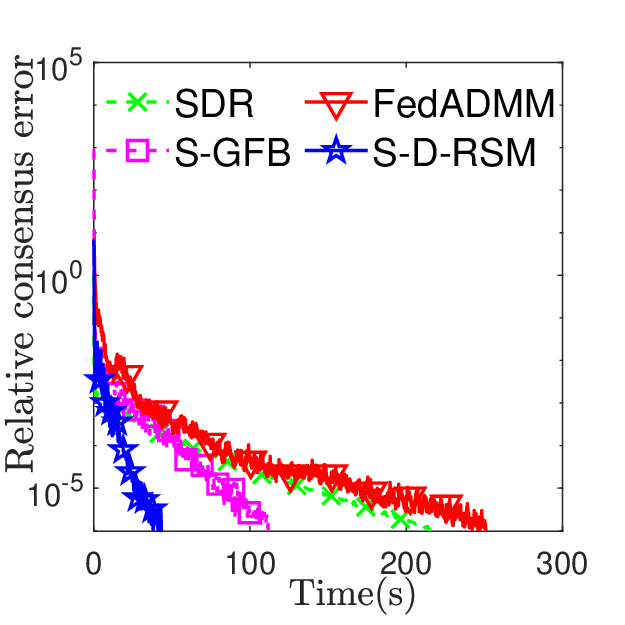} 
\end{minipage}
\begin{minipage}{0.49\columnwidth}
 \includegraphics[scale=0.42]{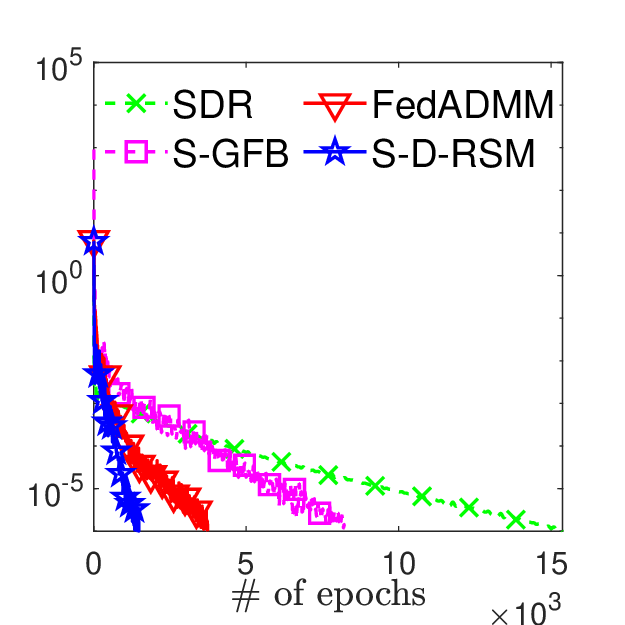} 
 \end{minipage}
\caption{Comparison of different methods for compressed sensing problems: DCT (top) and DFT (bottom).}
\label{fig1}
\end{figure}

We set \( n = 2500 \) and \( p = 0.25n \), and construct the sensing matrix \( A \) using the discrete cosine transform (DCT) or discrete Fourier transform (DFT). The ground-truth signal \( x^* \in \mathbb{R}^{2500} \) is generated from the standard normal distribution, with a sparsity level of 1\%. The observed vector is then computed as \( b = A x^* \). Figure \ref{fig1} shows the   relative consensus error $\max_{i}\left\{{\|y_i^k-x^k\|}/{\|x^k\|}\right\}$  of all the algorithms for the two sensing matrices $A$ (DCT,~DFT).
These figures clearly indicate that S-D-RSM converges much faster than SDR and S-GFB in terms of the objective value and relative error under all settings. Compared with S-GFB and FedADMM, the regularization term in S-D-RSM helps reduce the consensus error during the iteration process. Compared with SDR, the parallel update mechanism in S-D-RSM enhances computational efficiency, while the randomized selection of constraint sets \(C_i\) for projection reduces the frequency of redundant constraint processing in the linear system \(Ax = b\). 
\subsection{Logistic Regression Problem with \(\ell_1\)-norm Regularization Terms}
We further evaluate the performance of the proposed S-D-RSM algorithm on a logistic regression problem with \(\ell_1\)-norm regularization terms:
\[
\min_{x \in \mathbb{R}^n} \sum_{i=1}^{m} \frac{1}{m} \left( \log(1 + \exp(-b_i a_i^\top x)) + \lambda_i \|x\|_1 \right),
\]
where \(\{a_i\}_{i=1}^m \subset \mathbb{R}^n\) and \(\{b_i\}_{i=1}^m \subset \{\pm1\}\) denote the input features and output labels, respectively.

\begin{figure}[t]
\centering
\begin{minipage}{0.49\columnwidth}
 \includegraphics[scale=0.42]{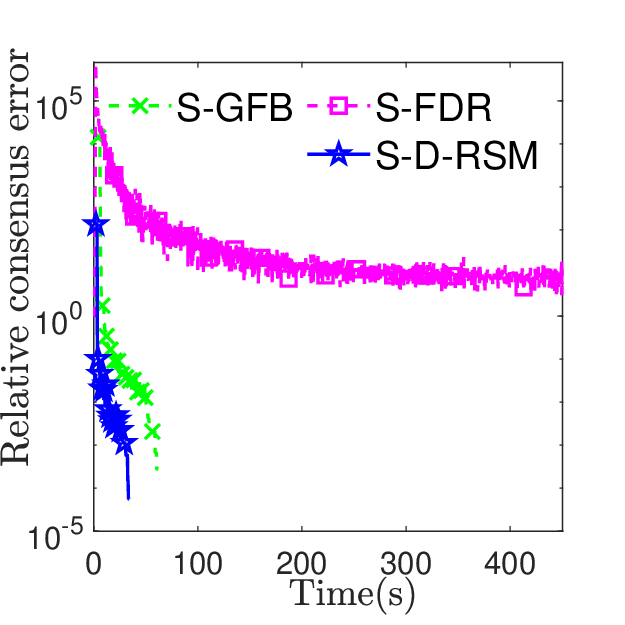} 
\end{minipage}
\begin{minipage}{0.49\columnwidth}
 \includegraphics[scale=0.42]{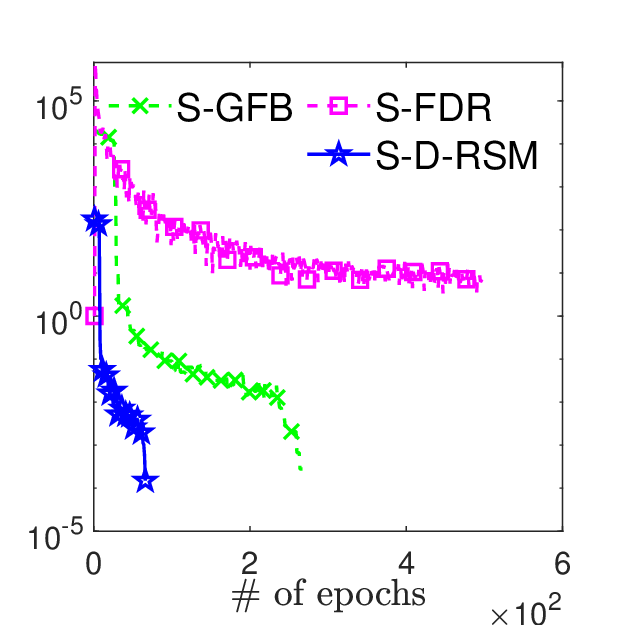} 
 \end{minipage}\\
 \begin{minipage}{0.49\columnwidth}
 \includegraphics[scale=0.42]{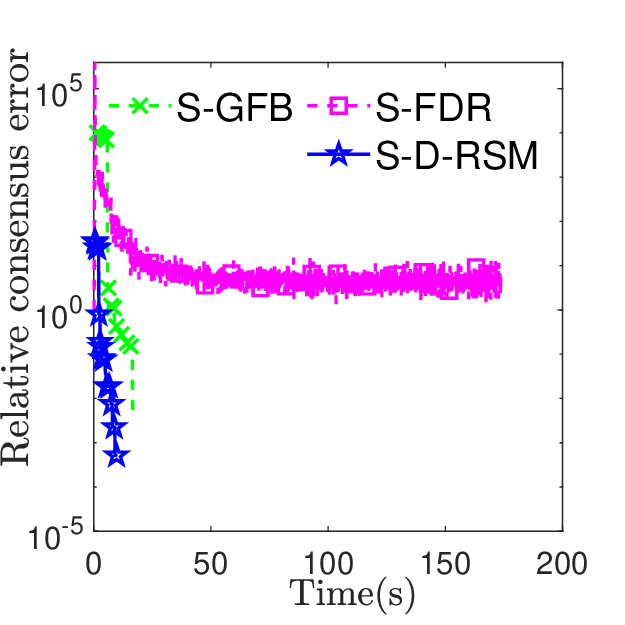} 
\end{minipage}
\begin{minipage}{0.49\columnwidth}
 \includegraphics[scale=0.42]{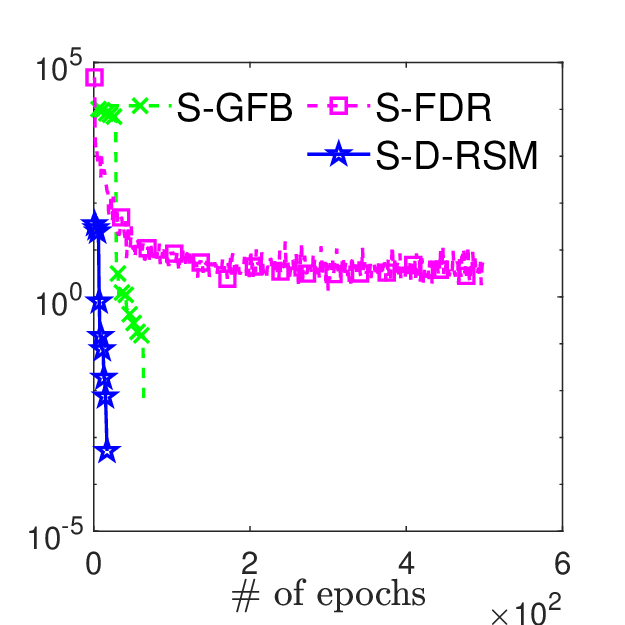} 
 \end{minipage}
\caption{Comparison of different methods for the logistic regression problem with $\ell_1$-norm regularizer on the two data sets: {a7a} (top) and  {mushrooms} (bottom).}
\label{fig3}
\end{figure}

Due to space limitations and the similarity of experimental patterns, we present results only for two commonly used benchmark datasets, namely {a7a} and {mushrooms}, obtained from the LIBSVM repository \cite{CC01a}. The datasets are randomly partitioned into 75\% training and 25\% testing sets and the maximum number of iterations is set to 1000. The regularization parameters \(\lambda_i\) are sampled uniformly from the interval \([10^{-3}, 10^{-2}]\).
S-D-RSM demonstrates a clear advantage over competing methods in terms of CPU time, as shown in Figure~\ref{fig3}. In particular, S-FDR employs a vanishing step size, which leads to slower convergence compared to S-GFB and S-D-RSM, both of which use constant step sizes. Although both S-GFB and S-D-RSM involve a comparable number of proximal mapping evaluations per iteration, S-D-RSM updates at least 30\% of gradients, whereas S-GFB computes the full gradient in every iteration. Consequently, S-D-RSM can potentially reduce the gradient computation cost by up to 70\% per iteration relative to S-GFB.
\section{Conclusions}
In this work, we propose a novel stochastic splitting algorithm, S-D-RSM, by integrating consensus optimization with operator splitting techniques. The method enables partial agent participation via   parallel updates and incorporates regularization to reduce consensus errors. In contrast to conventional stochastic methods, S-D-RSM is theoretically shown to achieve global convergence and an  $\mathcal{O}(\epsilon^{-1})$ complexity for both the objective value and consensus error, under constant step sizes and without strong convexity.
\section*{Acknowledgments}
The authors are grateful to the anonymous referees for their valuable comments and suggestions. This work was supported by the National Natural Science Foundation of China under Grant No.~12471290, and the Postgraduate Research \& Practice Innovation Program of Jiangsu Province under Grant No.~KYCX25\_1928.
\bibliography{sn-bib}
\onecolumn
\section{Appendices}
This appendix provides the detailed theoretical analysis for Algorithm~\ref{A}, referred to as the \textit{Stochastic Distributed Regularized Splitting Method (S-D-RSM)}. We begin by introducing some fundamental results that form the basis for the convergence and complexity analysis. Subsequently, we present complete proofs of all theoretical results related to S-D-RSM.
\paragraph{Problem Statement and Preliminary Results.}
In this work, we study the problem  large-scale distributed composite convex optimization:
\begin{equation}\label{p1}
  \min_{x\in \mathbb{R}^n} \left\{\Phi(x):=\sum_{i=1}^{m}\left(f_i(x)+g_i(x)\right)\right\},
\end{equation}
where $m$ is the numbers of nodes, $\{f_i\}_{i=1}^m$ is a sequence of proper, lower semicontinuous convex functions (not necessarily differentiable), $\{g_i\}_{i=1}^m$ is a sequence of convex functions that are Fr\'{e}chet differentiable on $\mathbb{R}^n$, and each gradient $\nabla g_i$ is $\frac{1}{\beta_i}$-Lipschitz continuous.

For all $x, y, z, w \in \mathbb{R}^n$, the following identity holds:
\begin{equation}\label{fp}
2\langle x - y, z - w \rangle = \|x - w\|^2 + \|y - z\|^2 - \|x - z\|^2 - \|y - w\|^2.
\end{equation}
We denote by $\Gamma_0(\mathbb{R}^n)$ the set of all proper, lower semicontinuous, and convex functions on $\mathbb{R}^n$. 
Given $f \in \Gamma_0(\mathbb{R}^n)$, the subdifferential of $f$ is defined as
\[
\partial f : x \mapsto \left\{ u \in \mathbb{R}^n \mid f(y) \ge f(x) + \langle u, y - x \rangle,\ \forall y \in \mathbb{R}^n \right\},
\]
and its proximal mapping is defined by
\[
\operatorname{prox}_f : x \mapsto \arg\min_{y \in \mathbb{R}^n} \left\{ f(y) + \frac{1}{2} \| y - x \|^2 \right\}.
\]
Using the definition of the  proximal mapping, we can immediately conclude that
\begin{equation}\label{dp}
y = \operatorname{prox}_f(x) ~ \Leftrightarrow ~ x - y \in \partial f(y).
\end{equation}  
Moreover, $ \operatorname{prox}_f$ is firmly nonexpansive (FNE)\cite[Theorem 6.42]{First2017}, i.e.,
\[
\langle x - z, \operatorname{prox}_f(x)- \operatorname{prox}_f(z) \rangle \geq \| \operatorname{prox}_f(x)- \operatorname{prox}_f(z) \|^2, \quad \forall x, z \in \mathbb{R}^n.
\]

The following lemma gives the properties of convex functions with Lipschitz continuous gradients.
\begin{Lemma}\label{le1}
Let $g \in \Gamma_0(\mathbb{R}^n)$ be Fr\'{e}chet differentiable with $L$-Lipschitz continuous gradient $\nabla g$. Then for all $x,y,z \in \mathbb{R}^n$, the following inequalities hold:
\begin{equation}\label{co}
\langle \nabla g(x) - \nabla g(y), x - z \rangle \geq -\frac{L}{4}\|y - z\|^2,
\end{equation}
and 
\begin{equation}\label{dc}
g(y) \leq g(z) + \langle \nabla g(x), y - z \rangle + \frac{L}{2}\|y - x\|^2.
\end{equation}
\end{Lemma}

\begin{proof}
\rm {Proof of (\ref{co}):} Using the Baillon-Haddad theorem \cite[Theorem 18.15]{Combettes2017Convex} and Young's inequality:
\begin{align}
\langle \nabla g(x) - \nabla g(y), x - z \rangle 
&= \langle \nabla g(x) - \nabla g(y), x - y \rangle + \langle \nabla g(x) - \nabla g(y), y - z \rangle \nonumber \\
&\geq \frac{1}{L}\|\nabla g(x) - \nabla g(y)\|^2 - \frac{1}{2\varepsilon}\|\nabla g(x) - \nabla g(y)\|^2 - \frac{\varepsilon}{2}\|y - z\|^2
\end{align}
for any $\varepsilon > 0$. Taking $\varepsilon = L/2$ yields the desired inequality.

{Proof of (\ref{dc}):} By the descent lemma \cite[Theorem 18.15]{Combettes2017Convex} and convexity of $g$:
\begin{align*}
g(y) - g(z) &= g(y) - g(x) + g(x) - g(z) \\
&\leq \langle \nabla g(x), y - x \rangle + \frac{L}{2}\|y - x\|^2 + \langle \nabla g(x), x - z \rangle \\
&= \langle \nabla g(x), y - z \rangle + \frac{L}{2}\|y - x\|^2.
\end{align*}
\end{proof}
The following two lemmas establish the convergence properties of a quasi-Fej\`{e}r monotone sequence of random variables. 
\begin{Lemma}{\rm \cite[Theorem 1]{ROBBINS1971233}}\label{s-m}
Let \( \{ \mathcal{F}_k \}_{k=0}^\infty \) be a sequence of \(\sigma\)-algebras such that \( (\forall k \in \mathbb{N}) \) \( \mathcal{F}_k \subset \mathcal{F}_{k+1} \). For each \( k \in \mathbb{N} \), let \( \chi_k, \beta_k, \gamma_k, \delta_k \) be nonnegative \( \mathcal{F}_k \)-measurable random variables such that\[
\mathbb{E}(\chi_{k+1} \mid \mathcal{F}_k) + \beta_k \le (1 + \gamma_k) \chi_k + \delta_k, \quad \mathbb{P}\text{-{\rm a.s.}},
\]
with
\[
\sum_{k=1}^\infty \gamma_k < \infty, \quad \mathbb{P}\text{-{\rm a.s.}}, \quad \sum_{k=1}^\infty  \delta_k < \infty, \quad \mathbb{P}\text{-{\rm a.s.}}.
\]
Then
\[
\sum_{k=1}^\infty  \beta_k < \infty, \quad \mathbb{P}\text{-{\rm a.s.}},
\]
and the sequence \( \{ \chi_k \}_{k \in \mathbb{N}} \) converges \(\mathbb{P}\)-almost surely to a random variable taking values in \( [0, +\infty) \).
\end{Lemma}
\begin{Lemma}{\rm{\cite[Proposition 2.3]{2015STOCHASTIC}}}\label{s-o}
Let \( F \subseteq \mathbb{R}^N \) be a nonempty closed set and \( \{ u^k \}_{k \in \mathbb{N}} \) an \( \mathbb{R}^N \)-valued random sequence. If for all \( u \in F \), there exists \( \bar{\Omega}_u \in \mathcal{F} \) with \( \mathbb{P}(\bar{\Omega}_u) = 1 \) such that \( \| u^k(\omega) - u \| \) converges for all \( \omega \in \bar{\Omega}_u \), then:

{\rm (i)}. There exists \( \hat{\Omega} \in \mathcal{F} \) such that \( \mathbb{P}(\hat{\Omega}) = 1 \), and for all \( \omega \in \hat{\Omega} \) and \( u \in F \), the sequence \( \{ \| u^k(\omega) - u \| \}_{k=0}^\infty \) converges.
    
{\rm (ii)}. If there exists \( \tilde{\Omega} \in \mathcal{F} \) such that \( \mathbb{P}(\tilde{\Omega}) = 1 \), and for every \( \omega \in \tilde{\Omega} \), every weak sequential cluster point of \( \{ u^k(\omega) \}_{k=1}^\infty \) belongs to \( F \), then the sequence \( \{ u^k \}_{k=1}^\infty \) converges weakly \( \mathbb{P} \)-{\rm a.s.} to an \( F \)-valued random variable.
\end{Lemma}

The following lemma shows that the outer semicontinuity of the subdifferential mapping of a proper, lower semicontinuous convex function.
\begin{Lemma}{\rm{\cite[Corollary 26.8]{Combettes2017Convex}}}\label{s-o1}
Let \( p \geq 2 \) be an integer and \( \{ G_i \}_{i=1}^p \subset \Gamma_0(\mathbb{R}^n) \). For each \( i \in [p] \), consider a sequence \( \{ (x_i^k, u_i^k) \}_{k=0}^\infty \subset \operatorname{gph} \partial G_i \) with \( u_i^k \in \partial G_i(x_i^k) \), and limit points \( (x_i, u_i) \in \mathbb{R}^n \times \mathbb{R}^n \). Suppose the following conditions hold:
\[
\sum_{i=1}^{p} u_i^k \to 0 \quad \text{and} \quad \forall i \in [p], \quad 
\begin{cases}
x_i^k \to x_i, \\
u_i^k \to u_i, \\
\sum_{t=1}^{p} (x_i^k - x_t^k) \to 0.
\end{cases}
\]
Then there exists \( x \in \operatorname{zer} \left( \sum_{i=1}^p \partial G_i \right) \) satisfying:

    {\rm (i)}. \( x = x_1 = \cdots = x_p \).
    
    {\rm (ii)}. \( \sum_{i=1}^{p} u_i = 0 \).
    
    {\rm (iii)}. For all \( i \in [p] \), $u_i\in\partial G_i(x)$.
\end{Lemma}

\paragraph{Stochastic Distributed Regularized Splitting Method and Its Theoretical Properties.}
Stochastic distributed regularized splitting method (S-D-RSM) can be described as follows.
\begin{algorithm}[H]
\caption{S-D-RSM for solving problem (\ref{p1})}
\label{A}
\textbf{Input}: $K>0$; $\alpha_i\ge0$; $\sigma\in[0,1]$; $\alpha_i+\lfloor1-\sigma\rfloor\neq0$; initial point $ y_i^0, z_i^0\in\mathbb{R}^n$, $i\in[m-1]$. \\
\textbf{Parameter}:  \hspace{-2pt} $\gamma\hspace{-2pt}\in\hspace{-2pt}\left(0,\min\limits_{i\in[m-1]}\left\{\frac{2\alpha_i}{\frac{1}{\beta_m(m-1)}+\frac{\sigma}{\beta_i}},\frac{2(2+\alpha_i)\beta_i}{1-\sigma}\right\}\right)$; $\lambda_i\in\left(0,2+\alpha_i-\frac{(1-\sigma)\gamma}{2\beta_i}\right)$; error$=1$ and the error tolerance $\varepsilon>0$.\\
\textbf{Output}: Approximate solution $x^k$.\\
\textbf{Process}:
\begin{algorithmic}[1]
\STATE Let $k=0$, error $=1$.
\WHILE{error $>\varepsilon$ or $k\le K$}
    \STATE Server update
 { \begin{align}\label{x1}
   x^{k+1}=& \operatorname{prox}_{\frac{\gamma f_m}{m-1}}\Bigg(\frac{1}{m-1}\sum_{i=1}^{m-1}\Big(z_i^k + \alpha_i(y_i^k-x^{k+1})-\frac{\gamma}{m-1}\nabla g_m(y_i^k)\Big)\\ &-\frac{\sigma\gamma}{m-1}\sum_{i=1}^{m-1}\nabla g_i(y_i^k)\Bigg)\nonumber.
    \end{align}}
    
    \STATE Randomly select users $S_k \subseteq [m-1]$  such that $S_k\perp\perp\mathcal{F}_k$  with $\mathbb{P}(i\in S_k)=p_i>0$, $i\in[m-1]$.    
    \STATE For user $i \in S_k$, compute
  {  \begin{align}
 y_i^{k+1}&= \operatorname{prox}_{\gamma f_i}\Big(2x^{k+1}-z_i^{k} + \alpha_i(x^{k+1}-y_i^{k+1}) - (1-\sigma)\gamma\nabla g_i(x^{k+1})\Big), \label{y1}\\
    z_i^{k+1}&= z_i^{k} + \lambda_i\left(y_i^{k+1}-x^{k+1}\right). 
    \end{align}}
    
    \STATE For user $i \notin S_k$, set
  { \begin{align*}
  y_i^{k+1} &= y_i^{k}, \\
    z_i^{k+1} &= z_i^{k}.
    \end{align*}}
    \STATE Update error $\leftarrow \frac{\sum_{i=1}^{m-1}\|y_i^k-x^k\|^2}{\|x^k\|^2}$ and $k \leftarrow k+1$.
\ENDWHILE
\end{algorithmic}
\end{algorithm}
In order to facilitate the convergence analysis of Algorithm \ref{A}, we introduce the following auxiliary variables, which are not computed in practice:
\begin{numcases}{ }
  \tilde{y}_i^{k+1}=\operatorname{prox}_{\gamma f_i}\left(2x^{k+1}-z_i^{k}+\alpha_i(x^{k+1}-\tilde{y}_i^{k+1})-{(1-\sigma)\gamma}\nabla g_i(x^{k+1})\right), ~i\in[m-1],\label{al2} \\
   \tilde{z}_i^{k+1}=z_i^{k}+\lambda_i\left(\tilde{y}_{i}^{k+1}-x^{k+1}\right),~i\in[m-1]\label{all4}.
\end{numcases}
The following lemma reveals the relation between the virtual user variables \( (\tilde{y}_i^k, \tilde{z}_i^k) \) and the actual user variables \( (y_i^k, z_i^k) \), which is fundamental for the convergence analysis.
 \begin{Lemma}\label{ce}
{ Let \( \left\{ x^k, ( y_i^k,z_i^k)_{i=1}^{m-1} \right\}_{k=0}^\infty \) denote the sequence generated by Algorithm \ref{A}, and let \( \left\{ (\tilde{y}_i^k,\tilde{z}_i^k)_{i=1}^{m-1} \right\}_{k=1}^\infty \) represent the auxiliary  variables defined by equations \eqref{al2} and \eqref{all4}. Then, for any \( x, z \in \mathbb{R}^n \), the following holds:
\begin{equation}\label{e2}
\mathbb{E}\left( \| y_i^{k+1} - x \|^2 \Big| \mathcal{F}_k \right) = p_i \| \tilde{y}_i^{k+1} - x \|^2 + (1 - p_i) \| y_i^k - x \|^2,
\end{equation}
and
\begin{equation}\label{e1}
\mathbb{E}\left( \| z_i^{k+1} - z \|^2 \Big| \mathcal{F}_k \right) = p_i \| \tilde{z}_i^{k+1} - z \|^2 + (1 - p_i) \| z_i^k - z \|^2, \quad i \in [m-1] ,
\end{equation}
where \( \mathcal{F}_k = \sigma \left( \boldsymbol{y}^0, \boldsymbol{y}^1, \dots, \boldsymbol{y}^k, \boldsymbol{z}^0, \boldsymbol{z}^1, \dots, \boldsymbol{z}^k \right) \), \( \boldsymbol{y}^k = (y_i^k)_{i=1}^{m-1} \), and \( \boldsymbol{z}^k = (z_i^k)_{i=1}^{m-1} \).
}
\end{Lemma}
\begin{proof}
{\rm 
For each \( i \in [m-1]  \), define the random variable $\varepsilon_k^i$ taking values 0 and 1 as follows:
\[
\varepsilon_k^i =
\begin{cases}
1, & \text{if } i \in S_k, \\
0, & \text{otherwise}.
\end{cases}
\]
It follows that  \( \mathbb{P}([\varepsilon_k^i = 1]) = \mathbb{P}(i \in S_k) = p_i \), and
\begin{numcases}{}
y_i^{k+1} = y_i^k + \varepsilon_k^i(\tilde{y}_i^{k+1} - y_i^k), & \( i \in [m-1]  \), \label{al3} \\
z_i^{k+1} = z_i^k + \lambda_i \varepsilon_k^i(\tilde{y}_i^{k+1} - x^{k+1}), & \( i \in [m-1]  \). \label{al4}
\end{numcases}

For every \( k \in \mathbb{N} \), let \( \mathcal{E}_k = (\varepsilon_k^1, \varepsilon_k^2, \dots, \varepsilon_k^{m-1}) \). The events \( \left\{ [{\mathcal{E}}_k = {\bf\mathcal{E}}] \right\}_{\mathcal{E} \in \mathcal{D}} \) form an almost sure partition of \( \Omega \), that is,
\[
\Omega = \bigcup_{\mathcal{E} \in \mathcal{D}} [{\mathcal{E}}_k = {\bf\mathcal{E}}], \quad \mathbb{P}\text{-a.s.}, \quad \text{and} \quad [{\mathcal{E}}_k = {\mathcal{E}_i}] \cap [{\mathcal{E}}_k = {\mathcal{E}_j}] = \emptyset, \, \forall i \neq j,
\]
where \( \mathcal{D} = \{0, 1\}^{m-1} \setminus \{(0, 0, \dots, 0)\} \). Therefore, we have
\[
\sum_{\mathcal{E} \in \mathcal{D}} I_{[{\mathcal{E}}_k = {\bf\mathcal{E}}]} = 1.
\]
Since the selection of \( S_k \) is independent of \( \{(z_i^j, y_i^j)_{i=1}^{m-1}\}_{j=0}^k \), the random variable \( \varepsilon_k^i \) is independent of \( \mathcal{F}_k \) for every \( i \in [m-1]  \).

By the independence of  \( \varepsilon_k^i \) and \( \mathcal{F}_k \), and the continuity of the norm, we have
\[
\begin{aligned}
\mathbb{E}\left( \| y_i^{k+1} - x \|^2 \Big| \mathcal{F}_k \right) &= \mathbb{E}\left( \| y_i^k + \varepsilon_k^i (\tilde{y}_i^{k+1} - y_i^k) - x \|^2 \Big| \mathcal{F}_k \right) \\
&= \mathbb{E}\left( \| y_i^k + \varepsilon_k^i (\tilde{y}_i^{k+1} - y_i^k) - x \|^2 \sum_{\mathcal{E} \in \mathcal{D}} I_{[{\mathcal{E}}_k = {\bf\mathcal{E}}]} \Bigg| \mathcal{F}_k \right) \\
&= \mathbb{E}\left( \sum_{\mathcal{E} \in \mathcal{D}} \| y_i^k + \varepsilon_k^i (\tilde{y}_i^{k+1} - y_i^k) - x \|^2 I_{[{\mathcal{E}}_k = {\bf\mathcal{E}}]} \Bigg| \mathcal{F}_k \right) \\
&= \sum_{\mathcal{E} \in \mathcal{D}} \mathbb{E}\left( I_{[{\mathcal{E}}_k = {\bf\mathcal{E}}]} \big| \mathcal{F}_k \right) \| y_i^k + \varepsilon_k^i (\tilde{y}_i^{k+1} - y_i^k) - x \|^2 \\
&= \sum_{\mathcal{E} \in \mathcal{D}} \mathbb{P}([{\mathcal{E}}_k = {\bf\mathcal{E}}]) \| y_i^k + \varepsilon_k^i (\tilde{y}_i^{k+1} - y_i^k) - x \|^2 \\
&= \sum_{\mathcal{E} \in \mathcal{D}, \varepsilon_k^i = 1} \mathbb{P}([{\mathcal{E}}_k = {\bf\mathcal{E}}]) \| \tilde{y}_i^{k+1} - x \|^2 + \sum_{\mathcal{E} \in \mathcal{D}, \varepsilon_k^i = 0} \mathbb{P}([{\mathcal{E}}_k = {\bf\mathcal{E}}]) \| y_i^k - x \|^2 \\
&= p_i \|\tilde{y}_i^{k+1} - x \|^2 + (1 - p_i) \| y_i^k - x \|^2.
\end{aligned}
\]
Note that \( z_i^{k+1} = z_i^k + \varepsilon_k^i (\tilde{z}_i^{k+1} - z_i^k) \). Repeating the above argument yields \eqref{al4}, which completes the proof.
}
\end{proof}
\section*{$\bullet$ Global Convergence Analysis of S-D-RSM}
Based on Lemma \ref{le1}--\ref{ce}, we now present all theoretical proofs related to the convergence  of Algorithm~\ref{A} in our article, starting with the reformulation and characterization of the solution set of problem~\eqref{p1}.

\paragraph{The proof of Lemma 1 (Reformulation I).}{\rm 
{(i)}. If $x^\star$ minimizes $\Phi$, then by Assumption~1 we have
\[
x^\star \in \operatorname{zer}\left(\sum_{i=1}^{m} (\partial f_i + \nabla g_i)\right).
\]
This implies that there exist \( a_i \in \partial f_i (x^\star) \) for \( i = 1, 2, \dots, m-1 \), such that
\[
\sum_{i=1}^{m} \left(a_i + \nabla g_i( x^\star)\right) = 0.
\]
For each \( i = 1, 2, \dots, m-1 \), define
\[
z_i^\star := x^\star - (1 - \sigma)\gamma\nabla g_i(x^\star) - \gamma a_i.
\]
Then, using~\eqref{dp}, we have
\[
x^\star = \operatorname{prox}_{\gamma f_i} \left(2x^\star - z_i^\star - (1 - \sigma)\gamma\nabla g_i(x^\star) \right),~i\in[m-1].
\]
Furthermore, it follows that
\[
\frac{1}{m-1} \sum_{i=1}^{m-1} z_i^\star 
= x^\star + \frac{\gamma}{m-1} a_m 
+ \frac{\gamma}{m-1} \nabla g_m (x^\star) 
+ \frac{\sigma \gamma}{m-1} \sum_{i=1}^{m-1}\nabla g_i(x^\star).
\]
Using~\eqref{dp} again, we obtain
\[
x^\star = \operatorname{prox}_{\frac{\gamma}{m-1}  f_m} \left( 
\frac{1}{m-1} \sum_{i=1}^{m-1} z_i^\star 
- \frac{\gamma}{m-1} \nabla g_m (x^\star) 
- \frac{\sigma \gamma}{m-1} \sum_{i=1}^{m-1}\nabla g_i(x^\star) 
\right).
\]
This completes the proof of part~(i).

{(ii)}. If \( (z_1^\star, z_2^\star, \dots, z_{m-1}^\star, x^\star) \in \mathcal{S} \), then
\[
\begin{cases}
x^\star = \operatorname{prox}_{\frac{\gamma}{m-1}  f_m} \left( \frac{1}{m-1} \sum_{i=1}^{m-1} z_i^\star - \frac{\gamma}{m-1} \nabla g_m (x^\star) - \frac{\sigma \gamma}{m-1} \sum_{i=1}^{m-1}\nabla g_i(x^\star) \right), \\
x^\star = \operatorname{prox}_{\gamma f_i} \left( 2x^\star - z_i^\star- (1 - \sigma) \gamma\nabla g_i(x^\star) \right), \quad i \in [m-1] .
\end{cases}
\]
Therefore, there exist \( a_i \in \partial f_i(x^\star) \) (\( i \in [m] \)) such that
\[
\begin{cases}
x^\star + \frac{\gamma}{m-1} a_m = \frac{1}{m-1} \sum_{i=1}^{m-1} z_i^\star - \frac{\gamma}{m-1} \nabla g_m (x^\star) - \frac{\sigma \gamma}{m-1} \sum_{i=1}^{m-1}\nabla g_i(x^\star), \\
x^\star + \gamma a_i = 2x^\star - z_i^\star - (1 - \sigma) \gamma\nabla g_i(x^\star), \quad i \in [m-1] .
\end{cases}
\]
Summing these equations yields
\[
0 = \sum_{i=1}^{m} \left(a_i +\nabla g_i(x^\star)\right) \in \sum_{i=1}^{m} \left(\partial f_i(x^\star) +\nabla g_i(x^\star)\right).
\]
The proof is completed.}

\begin{Remark}
{\rm{By leveraging the continuity of proximal mappings, it follows that $\mathcal{S}$ is a non-empty closed set in \( \mathbb{R}^{mn} \).}}
\end{Remark}

Next, we establish the descent property of the sequence generated by Algorithm~\ref{A}, which forms the foundation for proving the convergence of the proposed method.

\paragraph{ The proof Lemma 2 (Decreasing properties).}{\rm 
For any given $(z_1^\star,z_2^\star,\cdots,z_{m-1}^\star,x^\star)\in\mathcal{S}$, using FNE of the $\operatorname{prox}_{\frac{\gamma}{m-1}  f_m}$ and $\operatorname{prox}_{\gamma f_i}$, \eqref{x1} and \eqref{al2}, we have
\begin{equation}\label{m1}
\begin{aligned}
 &\sum_{i=1}^{m-1}\left\langle  x^{k+1}- x^\star,  z_i^k-z_i^\star+\alpha_i( y_i^k-x^{k+1})-\frac{\gamma}{m-1}\nabla g_m(y_i^k)+\frac{\gamma}{m-1}\nabla g_m(x^\star)\right\rangle\\
 &+\sigma\gamma\sum_{i=1}^{m-1}\langle  x^{k+1}- x^\star,\nabla g_i(x^\star)-\nabla g_i(y_i^k)\rangle\ge (m-1)\| x^{k+1}- x^\star\|^2
\end{aligned}
\end{equation}
and
\begin{equation}\label{m2}
 \begin{aligned}
 &2\langle \tilde{y}_i^{k+1}- x^\star,  x^{k+1}- x^\star\rangle-\langle\tilde{y}_i^{k+1}- x^{*}, z_i^k- z_i^\star\rangle+\alpha_i\langle\tilde{y}_i^{k+1}- x^{*}, x^{k+1}-\tilde{y}_i^{k+1}\rangle\vspace{1ex}\\
&+{(1-\sigma)\gamma}\langle\tilde{y}_i^{k+1}- x^\star,\nabla g_i(x^\star)-\nabla g_i(x^{k+1})\rangle\ge \|\tilde{y}_i^{k+1}- x^\star\|^2,~i\in[m-1].
 \end{aligned}
 \end{equation}
Summing $i$ from 1 to $m-1$ in both sides of \eqref{m2} and summing this inequality and \eqref{m1}, we obtain that
\begin{equation}\label{m3}
\begin{aligned}
\sum_{i=1}^{m-1}\| x^{k+1}-\tilde{y}_i^{k+1}\|^2\le&\sum_{i=1}^{m-1}\langle z_i^k- z_i^\star, x^{k+1}-\tilde{y}_i^{k+1}\rangle+\sum_{i=1}^{m-1}\alpha_i\langle\tilde{y}_i^{k+1}- x^\star, x^{k+1}-\tilde{y}_i^{k+1}\rangle\vspace{1ex}\\
&+\sum_{i=1}^{m-1}\alpha_i\langle  x^{k+1}- x^\star, y_i^k- x^{k+1}\rangle\\
&+{(1-\sigma)\gamma}\sum_{i=1}^{m-1}\langle{\tilde{y}}_i^{k+1}- x^\star,\nabla g_i(x^\star)-\nabla g_i(x^{k+1})\rangle\vspace{1ex}\\
&+{\sigma\gamma}\sum_{i=1}^{m-1}\langle  x^{k+1}- x^\star, \nabla g_i(x^\star)-\nabla g_i(y_i^k)\rangle\\
&+\frac{\gamma}{m-1}\sum_{i=1}^{m-1}\langle  x^{k+1}- x^\star, \nabla g_m(x^\star)-\nabla g_m(y_i^k)\rangle.
\end{aligned}
\end{equation}
Taking $g=g_i,~L=1/\beta_i,~x= x^\star$,~ $y= x^{k+1}$ and $z= \tilde{y}_i^{k+1}$ in \eqref{co}, we have 
\begin{equation}\label{co1}
\langle{\tilde{y}}_i^{k+1}- x^\star,\nabla g_i(x^\star)-\nabla g_i({x}^{k+1})\rangle\le\frac{1}{4\beta_i}\| x^{k+1}-\tilde{y}_i^{k+1}\|^2,~i\in[m-1].
\end{equation}
Similarly, we have 
\begin{equation}\label{co2}
\langle  x^{k+1}- x^\star, \nabla g_i(x^\star)-\nabla g_i(y_i^k)\rangle\le\frac{1}{4\beta_i}\|x^{k+1}-y_i^{k}\|^2,~i\in[m-1],
\end{equation}
and 
\begin{equation}\label{coo2}
\langle  x^{k+1}- x^\star, \nabla g_m(x^\star)-\nabla g_m(y_i^k)\rangle\le\frac{1}{4\beta_m}\|x^{k+1}-y_i^{k}\|^2,~i\in[m-1].
\end{equation}
Using \eqref{fp}, we can obtain that
\begin{equation}\label{fq1}
\langle\tilde{y}_i^{k+1}- x^\star, x^{k+1}-\tilde{y}_i^{k+1}\rangle=\frac{1}{2}\left(\| x^{k+1}- x^\star\|^2-\| x^{k+1}-\tilde{y}_i^{k+1}\|^2-\|\tilde{y}_i^{k+1}-x^\star\|^2\right).
\end{equation}
and
\begin{equation}\label{fq2}
\langle  x^{k+1}- x^\star, y_i^k- x^{k+1}\rangle=\frac{1}{2}\left(\| y_i^{k}- x^\star\|^2-\| x^{k+1}- y_i^{k}\|^2-\| x^{k+1}- x^\star\|^2\right).
\end{equation}
Substituting \eqref{co1}, \eqref{co2}, \eqref{coo2}, \eqref{fq1} and \eqref{fq2} into \eqref{m3} implies 
\begin{equation}\label{m4}
\begin{aligned}
\sum_{i=1}^{m-1}\langle z_i^k- z^\star,{\tilde{y}}_i^{k+1}- x^{k+1} \rangle\le&\frac{1}{2}\sum_{i=1}^{m-1}\alpha_i\left(\|y_i^{k}- x^\star\|^2-\|\tilde{y}_i^{k+1}- x^\star\|^2\right)\\
&-\frac{1}{2}\sum_{i=1}^{m-1}\left(\alpha_i-\frac{\gamma}{2\beta_m(m-1)}-\frac{\sigma\gamma}{2\beta_i}\right)\| x^{k+1}- y_i^{k}\|^2\\
&-\frac{1}{2}\sum_{i=1}^{m-1}\left(2+\alpha_i-\frac{(1-\sigma)\gamma}{2\beta_i}\right)\| x^{k+1}-\tilde{y}_i^{k+1}\|^2.
\end{aligned}
\end{equation}
Because $\tilde{z}_i^{k+1}=z_i^k+\lambda_i(\tilde{y}_i^{k+1}-x_i^{k+1})~(i\in[1,m-1])$, we can get
\begin{equation}\label{c1}
\begin{aligned}
\langle z_i^k- z^\star,{\tilde{y}}_i^{k+1}- x^{k+1} \rangle=&\frac{1}{\lambda_i}\langle  z_i^k- z_i^\star,\tilde{z}_i^{k+1}-z_i^k\rangle\\
=&\frac{1}{2\lambda_i}\left(\|\tilde{z}_i^{k+1}- z_i^\star\|^2-\|z_i^{k}- z_i^\star\|^2\right)\hspace{-1.5pt}-\frac{\lambda_i}{2}\|{\tilde{y}}_i^{k+1}- x^{k+1}\|^2,i\in[m-1].
\end{aligned}
\end{equation}
Substituting \eqref{c1} into \eqref{m4} implies that 
\begin{equation}\label{m40}
\begin{aligned}
\sum_{i=1}^{m-1}\frac{1}{\lambda_i}\left(\|\tilde{z}_i^{k+1}- z_i^\star\|^2-\|z_i^{k}- z_i^\star\|^2\right)\le&\sum_{i=1}^{m-1}\alpha_i\left(\|y_i^{k}- x^\star\|^2-\|\tilde{y}_i^{k+1}- x^\star\|^2\right)\\
&-\sum_{i=1}^{m-1}\left(\alpha_i-\frac{\gamma}{2\beta_m(m-1)}-\frac{\sigma\gamma}{2\beta_i}\right)\| x^{k+1}- y_i^{k}\|^2\\
&-\sum_{i=1}^{m-1}\left(2+\alpha_i-\frac{(1-\sigma)\gamma}{2\beta_i}-\lambda_i\right)\| x^{k+1}-\tilde{y}_i^{k+1}\|^2.
\end{aligned}
\end{equation}
Combining \eqref{e2} and \eqref{e1}, we can get
\begin{equation}\label{e3}
\|\tilde{z}_i^{k+1}-z_i^\star\|^2-\|z_i^{k}-z_i^\star\|^2=\frac{1}{p_i}\mathbb{E}\left(\|z_i^{k+1}-z_i^\star\|^2 \mid\mathcal{F}_k\right)-\frac{1}{p_i}\|z_i^k-z_i^\star\|^2
\end{equation}
and
\begin{equation}\label{e4}
 \|y_i^{k}-x^\star\|^2- \|\tilde{y}_i^{k+1}-x^\star\|^2=\frac{1}{p_i}\|y_i^k-x^\star\|^2-\frac{1}{p_i}\mathbb{E}\left(\|y_i^{k+1}-x^\star\|^2\mid \mathcal{F}_k\right).
\end{equation}
Substituting  \eqref{e3} and \eqref{e4} into \eqref{m40}, we have 
\begin{equation}\label{m41}
\begin{aligned}
&\sum_{i=1}^{m-1}\left\{\frac{1}{\lambda_ip_i}\mathbb{E}\left(\|z_i^{k+1}-z_i^\star\|^2 \big|\mathcal{F}_k\right)-\frac{1}{\lambda_ip_i}\|z_i^k-z_i^\star\|^2\right\}\\
\le&\sum_{i=1}^{m-1}\left\{\frac{\alpha_i}{p_i}\|y_i^k-x^\star\|^2-\frac{\alpha_i}{p_i}\mathbb{E}\left(\|y_i^{k+1}-x^\star\|^2 \big|\mathcal{F}_k\right)\right\}\\
&-\sum_{i=1}^{m-1}\left(\alpha_i-\frac{\gamma}{2\beta_m(m-1)}-\frac{\sigma\gamma}{2\beta_i}\right)\| x^{k+1}- y_i^{k}\|^2\\
&-\sum_{i=1}^{m-1}\left(2+\alpha_i-\frac{(1-\sigma)\gamma}{2\beta_i}-\lambda_i\right)\| x^{k+1}-\tilde{y}_i^{k+1}\|^2.
\end{aligned}
\end{equation}
Reorganizing the terms conclude the proof.
}

\begin{Remark}
 {\rm  
 (i). In the case where $\nabla g_1=\nabla g_2=\cdots=\nabla g_m=0$, and applying a proof method analogous to the one previously discussed, we arrive at an inequality similar to \eqref{m40}:
 \begin{equation}\label{md0}
\begin{aligned}
& \mathbb{E}\left[ \sum_{i=1}^{m-1} \frac{1}{\lambda_i p_i} \| z_i^{k+1} - z_i^\star \|^2  \big| \mathcal{F}_k \right] + \mathbb{E}\left[ \sum_{i=1}^{m-1} \frac{\alpha_i}{p_i} \| y_i^{k+1} - x^\star \|^2  \big| \mathcal{F}_k \right] \\
\le & \sum_{i=1}^{m-1} \frac{1}{\lambda_i p_i} \| z_i^k - z_i^\star \|^2 + \sum_{i=1}^{m-1} \frac{\alpha_i}{p_i} \| y_i^k - x^\star \|^2 - \sum_{i=1}^{m-1}  \alpha_i  \| x^{k+1} - y_i^k \|^2 \\
& - \sum_{i=1}^{m-1} \left( 2 + \alpha_i - \lambda_i \right) \| x^{k+1} - \tilde{y}_i^{k+1} \|^2.
\end{aligned}
\end{equation}
 (ii). If \( S_k \equiv [m-1] \) for each \( k \in \mathbb{N} \), then \eqref{m40} can be rewritten as
  \begin{equation}\label{md}
\begin{aligned}
&\sum_{i=1}^{m-1}\left(\frac{1}{\lambda_i}\|z_i^{k+1}-z_i^\star\|^2+\alpha_i\|y_i^{k+1}-x^\star\|^2\right)\\
\le&\sum_{i=1}^{m-1}\left(\frac{1}{\lambda_i}\|z_i^{k+1}-z_i^\star\|^2+\alpha_i\|y_i^k-x^\star\|^2\right)\\
&-\sum_{i=1}^{m-1}\left(\alpha_i-\frac{\gamma}{2\beta_m(m-1)}-\frac{\sigma\gamma}{2\beta_i}\right)\| x^{k+1}- y_i^{k}\|^2\\
&-\sum_{i=1}^{m-1}\left(2+\alpha_i-\frac{(1-\sigma)\gamma}{2\beta_i}-\lambda_i\right)\| x^{k+1}-{y}_i^{k+1}\|^2.
\end{aligned}
\end{equation} }
\end{Remark}

Leveraging the problem reformulation in Lemma~1 and the descent property established in Lemma~2, we now present the proof of the global convergence of Algorithm~\ref{A}.

\paragraph{ The proof of Theorem 1 (Convergence).}{\rm {\rm (i).} For  any given $(z_1^\star, \cdots, z_{m-1}^\star, x^\star) \in \mathcal{S}$, define the nonnegative $\mathcal{F}_k$-measurable random variables as follows:
$$
\begin{cases}
\chi_k = \sum_{i=1}^{m-1} \left( \frac{1}{\lambda_i p_i} \| z_i^k - z_i^\star \|^2 + \frac{\alpha_i}{p_i} \| y_i^k - x^\star \|^2 \right),\vspace{1ex} \\
\beta_k = \sum_{i=1}^{m-1} \left( \alpha_i - \frac{\gamma}{2\beta_m(m-1)} - \frac{\sigma \gamma}{2\beta_i} \right) \| x^{k+1} - y_i^k \|^2 \\
~~~~~~~+ \sum_{i=1}^{m-1} \left( 2 + \alpha_i - \frac{(1-\sigma) \gamma}{2\beta_i} - \lambda_i \right) \| x^{k+1} - \tilde{y}_i^{k+1} \|^2,~\forall k\in\mathbb{N}.
\end{cases}
$$
From Lemma \ref{s-m} and Lemma 2, we deduce that the sequences $\{\chi_k\}_{k=0}^\infty$ and $ \sum_{k=1}^{\infty} \beta_k$ converge $\mathbb{P}$-a.s.  Hence, (i) follows.

{\rm (ii).} For each $k \in \mathbb{N}$, define nonnegative $\mathcal{F}_k$-measurable random variable by
$$
\chi'_k = \sum_{i=1}^{m-1} \left( \frac{1}{\lambda_i p_i} \| z_i^k - z_i^\star \|^2 \right) + c \| x^{k+1} - x^\star \|^2,
$$
where $c = \sum_{i=1}^{m-1} \frac{\alpha_i}{p_i}$ is a constant. For any given \( (z_1^\star, z_2^\star, \dots, z_{m-1}^\star, x^\star) \in \mathcal{S} \), by invoking (i), we conclude that the sequence \( \{\chi'_k\}_{k=0}^\infty \) converges  $\mathbb{P}$-a.s. 

Define an inner product $\ll \cdot, \cdot \gg$ on $\mathbb{R}^{mn}$ by
$$
(\forall \boldsymbol{x}, \boldsymbol{y} \in \mathbb{R}^{mn}) \quad \ll \boldsymbol{x}, \boldsymbol{y} \gg = \sum_{i=1}^{m-1} \frac{1}{\lambda_i p_i} \langle x_i, y_i \rangle + c \langle x_m, y_m \rangle.
$$
Denote the associated norm by $\mid \cdot \mid = \sqrt{\ll \cdot, \cdot \gg}$. Let $\mathbf{w}^k = (z_1^k, z_2^k, \cdots, z_{m-1}^k, x^{k+1})$ and $\mathbf{w}^\star = (z_1^\star, z_2^\star, \cdots, z_{m-1}^\star, x^\star) \in \mathcal{S}$, we have $\chi_k' = \mid \mathbf{w}^k - \mathbf{w}^\star \mid^2$. By utilizing the conclusion (i) of Lemma \ref{s-o} ($N=mn$), we can conclude that there exists $\hat{\Omega} \in \mathcal{F}$ such that $\mathbb{P}(\hat{\Omega}) = 1$ and, for all $\omega \in \hat{\Omega}$ and $\mathbf{w}^\star\in\mathcal{S}$, $\{\mid \mathbf{w}^k(\omega) - \mathbf{w}^\star \mid \}_{k=1}^\infty$ converges. Combining (i), we can know that there exists $\bar{\Omega}\in\mathcal{F}$ such that $\mathbb{P}(\bar{\Omega})=1$, and for all $\omega\in \bar{\Omega}$, 
\begin{equation}\label{00}
\lim\limits_{k \to \infty} \left\| x^{k+1}(\omega) - y_i^k(\omega) \right\| = \lim\limits_{k \to \infty} \| x^k(\omega) - \tilde{y}_i^k(\omega) \| = 0,~i \in [1,m-1].
\end{equation}
Let \( \tilde{\Omega} = \bar{\Omega} \cap \hat{\Omega} \). Suppose \( \tilde{\Omega} = \emptyset \). Then, by the additivity of probability measures for disjoint events, we would have  
\[
1 \geq \mathbb{P}(\bar{\Omega} \cup \hat{\Omega}) = \mathbb{P}(\bar{\Omega}) + \mathbb{P}(\hat{\Omega}) = 2 > 1,
\]  
which leads to a contradiction. Hence, \( \tilde{\Omega} \) must be nonempty, and it follows that  
\[
\mathbb{P}(\tilde{\Omega}) = \mathbb{P}\left(\bar{\Omega} -(\bar{\Omega} \cap \hat{\Omega}^c)\right)=\mathbb{P}(\bar{\Omega}) -\mathbb{P}(\bar{\Omega} \cap \hat{\Omega}^c)=1-0=1 .
\]  
Consequently, for all \( \omega \in \tilde{\Omega} \) and \( \mathbf{w}^\star \in \mathcal{S} \), the sequence \( \left\{ |\mathbf{w}^k(\omega) - \mathbf{w}^\star| \right\}_{k=1}^\infty \) converges, and \eqref{00} holds.  

Furthermore, the following boundedness properties hold for all \( \omega \in \tilde{\Omega} \):  \\
    (1).  The sequences \( \{ z_i^k(\omega) \}_{k=1}^\infty \), \( \{ \tilde{y}_i^k(\omega) \}_{k=1}^\infty \), \( \{{y}_i^k(\omega) \}_{k=1}^\infty \), and \( \{ x^{k+1}(\omega) \}_{k=1}^\infty \) are bounded in \( \mathbb{R}^n \) for each \( i \in [m-1] \). \\ 
    (2). Since each operator \( \nabla g_i \) (\( i \in [m-1] \)) is  continuous, the sequences \( \{\nabla g_i\left(x^{k+1}(\omega)\right) \}_{k=1}^\infty \) and \( \{ \nabla g_i \left(y_i^k(\omega)\right) \}_{k=1}^\infty \) are bounded in \( \mathbb{R}^n \).  \\
    (3). Similarly, \( \{ \nabla g_m \left(y_i^k(\omega)\right) \}_{k=1}^\infty \) remains bounded in \( \mathbb{R}^n \) for all \( i \in [m-1] \)(due to continuity).

For any given $\omega \in \tilde{\Omega}$, let $(z_1^\infty(\omega), z_2^\infty(\omega), \cdots, z_{m-1}^\infty(\omega), x^\infty(\omega)) \in \mathbb{R}^{mn} $ be a sequential cluster point of $\{\mathbf{w}^k(\omega)\}$. Without loss of generality, for each $i \in [m-1]$ and $\omega \in \tilde{\Omega}$, we may assume that
$$
z_i^{k_j}(\omega) \to z_i^\infty(\omega), \quad x^{k_j+1}(\omega) \to x^\infty(\omega), \quad j \to \infty.
$$
For each $i \in [m-1]$ and $\omega \in \tilde{\Omega}$, it is clear from conclusion (i) that
\begin{equation}\label{i0}
\tilde{y}_i^{k_j+1}(\omega) \to x^\infty(\omega), \quad j \to \infty.
\end{equation}
For all $i\in[m-1]$ and $j\in\mathbb{N}$, by utilizing the definition of the resolvent, \eqref{x1} and \eqref{y1}, we have  
\begin{equation}\label{i1}
\begin{cases}
    u_m^{j}:=\sum_{i=1}^{m-1}\Bigg\{z_i^{k_j}+\alpha_i\left(y_i^{k_j}-x^{k_j+1}\right)-x^{k_j+1}+\frac{\gamma}{m-1}\left(\nabla g_m\left(x^{k_j+1}\right)-\nabla g_m\left(y_i^{k_j}\right)\right)\\
   ~~~~~~~ +\sigma\gamma\left(\nabla g_i\left(x^{k_j+1}\right)-\nabla g_i\left(y_i^{k_j}\right)\right)\Bigg\}\\
    ~~~~~\in \gamma\left(\partial f_m+\sigma \sum_{i=1}^{m-1}\nabla g_i+\nabla g_m\right)\left(x^{k_j+1}\right),\vspace{1ex} \\
    u_i^{j}:=2x^{k_j+1}-z_i^{k_j}+\alpha_i\left(x^{k_j+1}-\tilde{y}_i^{k_j+1}\right)-\tilde{y}^{k_j+1}_i+{(1-\sigma)\gamma}\left(\nabla g_i\left(\tilde{y}^{k_j+1}_i\right)-\nabla g_i\left(x^{k_j+1}\right)\right)\\
    ~~~~~\in\gamma\left(\partial f_i+(1-\sigma)\nabla g_i\right)\left(\tilde{y}_i^{k_j+1}\right).
  \end{cases}
\end{equation}
Hence we have $$\sum_{i=1}^{m}u_i^{j}=\sum_{i=1}^{m-1}\left\{x^{k_j+1}-\tilde{y}_i^{k_j+1}+\alpha_i\left(y_i^{k_j}-\tilde{y}_i^{k_j+1}\right)+\Delta_i^{j}+\Psi_i^{j}\right\},$$
where $\Delta_i^{j}={\sigma\gamma}\left(\nabla g_i\left(x^{k_j+1}\right)-\nabla g_i\left(y_i^{k_j}\right)\right)+\frac{\gamma}{m-1}\left(\nabla g_m\left(x^{k_j+1}\right)-\nabla g_m\left(y_i^{k_j}\right)\right)$ and $\Psi_i^{j}={(1-\sigma)\gamma}\left(\nabla g_i\left(\tilde{y}^{k_j+1}_i\right)-\nabla g_i\left(x^{k_j+1}\right)\right).$ Because $$\lim\limits_{j\to\infty}\left\|x^{k_j+1}(\omega)-y_i^{k_j}(\omega)\right\|=\lim\limits_{j\to\infty}\left\|x^{k_j+1}(\omega)-\tilde{y}_i^{k_j+1}(\omega)\right\|=0~\left(\forall \omega\in\tilde{\Omega},i\in[m-1]\right)$$ and $\nabla g_i$ is $\beta_i^{-1}$ Lipschitz continuous, we have $$\lim\limits_{j\to\infty}\Delta_i^{j}(\omega)=\lim\limits_{j\to\infty}\Psi_i^{j}(\omega)=0~\left(\forall \omega\in\tilde{\Omega},i\in[m-1]\right)$$ and
\begin{equation}\label{i3}\begin{cases}
\lim\limits_{j\to\infty}\left\|y^{k_j}(\omega)-\tilde{y}_i^{k_j+1}(\omega)\right\|=0,\vspace{1ex} \\
u_m^j(\omega)\to \sum_{i=1}^{m-1} \left[z_i^\infty(\omega)-x^\infty(\omega)\right]:=u_m^\infty(\omega),\forall \omega\in\tilde{\Omega},\vspace{1ex} \\
u_i^j(\omega)\to x^\infty(\omega)-z_i^\infty(\omega):=u_i^\infty(\omega), \forall \omega\in\tilde{\Omega},~i\in[m-1].
\end{cases}
\end{equation}
For all $\omega\in\tilde{\Omega}$, it holds that 
$$\left\|\sum_{i=1}^{m}u_i^{j}(\omega)\right\|\le \sum_{i=1}^{m}\left(\left\|x^{k_j+1}(\omega)-\tilde{y}_i^{k_j+1}(\omega)\right\|+\alpha_i\left\|y_i^{k_j}(\omega)-\tilde{y}_i^{k_j+1}(\omega)\right\|+\left\|\Delta_i^{j}(\omega)\right\|+\left\|\Phi_i^{j}(\omega)\right\|\right)\to 0,$$that is \begin{equation}\label{i2}
  \lim_{j\to\infty}\sum_{i=1}^{m}u_i^j(\omega)=0,~\forall \omega\in\tilde{\Omega}.
\end{equation}
Denote $\tilde{y}_m^k=x^k(\forall k\in\mathbb{N})$, thus, for all $\omega\in\tilde{\Omega}$ and $i\in[m]$, it holds that 
\begin{equation}\label{i4}
 \lim\limits_{j\to\infty}\sum_{t=1}^{m}\left(\tilde{y}_i^{k_j+1}(\omega)-\tilde{y}_t^{k_j+1}(\omega)\right)=0
 \end{equation} 
from \eqref{00}. Combing \eqref{i0}-\eqref{i4} and using  Lemma \ref{s-o1},  we can know that for any given $\omega\in\tilde{\Omega}$, $x^\infty(\omega)\in$ zer$\left(\sum_{i=1}^{m}(\partial f_i+\nabla g_i)\right)$ such that 
\begin{equation}\label{cl}
\begin{cases}
\left(x^\infty(\omega),u_i^\infty(\omega)\right)\in {\rm{gra}}\left(\gamma\left( \partial f_i+(1-\sigma)\nabla g_i\right)\right),~\forall i\in[m-1],\vspace{1ex}\\
\left(x^\infty(\omega),u_m^\infty(\omega)\right)\in {\rm{gra}}\left(\gamma\left( \partial f_m+\sigma \sum_{i=1}^{m-1}\nabla g_i+\nabla g_m\right)\right).
\end{cases}
\end{equation}
Combining \eqref{i3}, we can rewrite \eqref{cl} as 
\begin{equation}\label{cl1}
\begin{cases}
 x^\infty(\omega)-z_i^\infty(\omega)\in \gamma\left( \partial f_i+(1-\sigma)\nabla g_i\right)\left(x^\infty(\omega)\right), \forall \omega\in\tilde{\Omega},~i\in[m-1],\vspace{1ex} \\
\sum_{i=1}^{m-1} \left(z_i^\infty(\omega)-x^\infty(\omega)\right)\in \gamma\left( \partial f_m+\sigma \sum_{i=1}^{m-1}\nabla g_i+\nabla g_m\right)\left(x^\infty(\omega)\right),\forall \omega\in\tilde{\Omega}.
\end{cases}
\end{equation}
Utilizing the definition of resolvent, the inclusion \eqref{cl1} can be reformulated as
$$\begin{cases}
x^\infty(\omega)=\operatorname{prox}_{\gamma f_i}\left(2x^\infty(\omega)- z_i^\infty(\omega)-{(1-\sigma)\gamma} \nabla g_i(x^\infty(\omega))\right),~\forall \omega\in\tilde{\Omega}, i\in[m-1],\vspace{1ex}\\
x^\infty(\omega)=\operatorname{prox}_{\frac{\gamma f_m}{m-1}}\left(\frac{1}{m-1}\sum_{i=1}^{m-1}\left(z_i^\infty(\omega)-{\sigma\gamma}\nabla g_i(x^\infty(\omega))\right)-\frac{\gamma}{m-1} \nabla g_m(x^\infty(\omega))\right), \forall \omega\in\tilde{\Omega}. 
\end{cases}$$ 
Thus, for all $\omega\in\tilde{\Omega}$ where $\mathbb{P}(\tilde{\Omega})=1$, we can get $\left(z_1^\infty(\omega),z_2^\infty(\omega),\cdots,z_{m-1}^\infty(\omega),x^\infty(\omega)\right)\in\mathcal{S}$ from the definition of $\mathcal{S}$. Using Lemma \ref{s-o}, we can know that there exists $\mathcal{S}-$valued random variable $(\tilde{z}_1,\tilde{z}_2,\cdots,\tilde{z}_{m-1},\tilde{x})$, such that $x^k \to  \tilde{x}$ and $z_i^k\to\tilde{z}_i~(i\in[m-1]),~\mathbb{P}$-a.s. Thus it holds that $x^k \to \tilde{x}\in$ zer$\left(\sum_{i=1}^{m}\left(\partial f_i+\nabla g_i\right)\right)$ $\mathbb{P}$-a.s. from Lemma 1. The proof is completed.}

\begin{Remark}
 {\rm If $S_k\equiv[m-1]$, by using inequality \eqref{md} and repeating the above proof process, we can obtain that under the same parameter conditions, there exists an \(x\in\mathrm{zer}\left(\sum_{i = 1}^{m}\left(\partial f_i + \nabla g_i\right)\right)\) such that \(x^k\to x\).}
\end{Remark}
\section*{$\bullet$ Complexity Analysis of S-D-RSM}
We now proceed to the complexity analysis of Algorithm~\ref{A}. As a preliminary step, we establish the relationship between the two problem reformulations-namely, the connection between the set \( \mathcal{S} \) in Lemma 1 and the original objective function.

\paragraph{The proof of Lemma 4.}{\rm From the definition of $\mathcal{S}$, it follows that
$$
\begin{cases}
\frac{x^\star-z_i^\star}{\gamma}-(1-\sigma)\nabla g_i(x^\star)\in \partial f_i(x^\star),~i\in[m-1],\\
\frac{1}{\gamma}\sum_{i=1}^{m-1}(z_i^\star-x^\star)-\nabla g_m(x^\star)-\sigma\sum_{i=1}^{m-1}\nabla g_i(x^\star)\in \partial f_m(x^\star).
\end{cases}
$$
Using the definition of the subdifferential of a convex function, we obtain
\begin{align}
f_i(x_i)\ge& f_i(x^\star)+\frac{1}{\gamma}\langle x^\star-z_i^\star,x_i -x^\star\rangle-(1-\sigma)\langle\nabla g_i(x^\star), x_i-x^\star\rangle,~ \forall x_i \in \mathbb{R}^n,~i\in[m-1], \label{si1}\\
f_m(x_m)\ge& f_m(x^\star)+\frac{1}{\gamma}\sum_{i=1}^{m-1}\langle z_i^\star-x^\star,x_m-x^\star  \rangle-\langle \nabla g_m(x^\star),x_m-x^\star \rangle \nonumber\\
&-\sigma\sum_{i=1}^{m-1}\langle \nabla g_i(x^\star),x_m-x^\star  \rangle,\quad \forall x_m \in \mathbb{R}^n. \label{si2}
\end{align}

Summing \eqref{si1} and \eqref{si2} yields
\begin{equation}\label{fc}
\begin{aligned}
\sum_{i=1}^{m} f_i(x_i)\ge& \sum_{i=1}^{m}f_i(x^\star)+\frac{1}{\gamma}\sum_{i=1}^{m-1}\langle z_i^\star-x^\star,x_m-x_i\rangle-(1-\sigma)\sum_{i=1}^{m-1}\langle\nabla g_i(x^\star), x_i-x^\star\rangle\\
& -\langle \nabla g_m(x^\star),x_m-x^\star \rangle-\sigma\sum_{i=1}^{m-1}\langle \nabla g_i(x^\star),x_m-x^\star  \rangle,~\forall x_i\in\mathbb{R}^n,~i\in[m].
\end{aligned}
\end{equation}

By the convexity of $g_i$, it follows that for each $i \in [m]$,
\begin{equation}\label{gc}
g_i(x_i)\ge g_i(x^\star) + \left\langle \nabla g_i(x^\star), x_i-x^\star \right\rangle,~ \forall x_i\in\mathbb{R}^n.
\end{equation}

Substituting \eqref{gc} into \eqref{fc}, we obtain
\begin{equation}\label{fc1}
\begin{aligned}
\sum_{i=1}^{m} f_i(x_i)\ge& \sum_{i=1}^{m}f_i(x^\star)+\sum_{i=1}^{m} g_i(x^\star)+\frac{1}{\gamma}\sum_{i=1}^{m-1}\langle z_i^\star-x^\star,x_m-x_i\rangle\\
& -\sum_{i=1}^{m-1}\left((1-\sigma)g_i(x_i)+\sigma g_i(x_m)\right)-g_m(x_m),~ \forall x_i\in\mathbb{R}^n,~i\in[m],
\end{aligned}
\end{equation}
which completes the proof.
 }

Next, we analyze the evolution of the objective function of problem \eqref{p1} along the sequence of iterates generated by Algorithm \eqref{A}.

\paragraph{The proof of Lemma 5.}{\rm Based on the first-order optimality condition for subproblems \eqref{x1} and \eqref{al2}, we have
\begin{equation*}
\sum_{i=1}^{m-1}\left(z_i^k-x^{k+1}+\alpha_i\left(y_i^k-x^{k+1}\right)-\frac{\gamma}{m-1}\nabla g_m(y_i^k)\right)-{\sigma\gamma}\sum_{i=1}^{m-1}\nabla g_i(y_i^k)\in {\gamma}\partial f_m(x^{k+1}),
\end{equation*}
and 
\begin{equation*}
2x^{k+1}-z_i^k+\alpha_i\left(x^{k+1}-\tilde{y}_i^{k+1}\right)-\tilde{y}_i^{k+1}-(1-\sigma)\gamma\nabla g_i(x^{k+1})\in\gamma \partial f_i(\tilde{y}_i^{k+1}),~i\in [m-1].
\end{equation*}
For any $(z_1^\star, z_2^\star, \dots, z_{m-1}^\star, x^\star) \in \mathcal{S}$, applying the definition of the subdifferential operator for a convex function, we conclude
\begin{equation*}
\begin{aligned}
f_m(x^\star)\ge& f_m(x^{k+1})+\frac{1}{\gamma}\sum_{i=1}^{m-1}\langle x^\star-x^{k+1}, z_i^k-x^{k+1} \rangle+\frac{1}{\gamma}\sum_{i=1}^{m-1}\alpha_i\langle x^\star-x^{k+1}, y_i^k-x^{k+1} \rangle\\
&-\frac{1}{m-1}\sum_{i=1}^{m-1}\langle x^\star-x^{k+1}, \nabla g_m(y_i^k) \rangle-\sigma\sum_{i=1}^{m-1}\langle x^\star-x^{k+1}, \nabla g_i(y_i^k) \rangle,
\end{aligned}
\end{equation*}
\begin{equation*}
\begin{aligned}
f_i(x^\star)\ge& f_i(\tilde{y}_i^{k+1})+\frac{1}{\gamma}\langle x^\star-\tilde{y}_i^{k+1},x^{k+1}-z_i^k \rangle+\frac{\alpha_i}{\gamma}\langle x^\star-\tilde{y}_i^{k+1}, x^{k+1}-\tilde{y}_i^{k+1} \rangle\\
&+\frac{1}{\gamma}\langle x^\star-\tilde{y}_i^{k+1},x^{k+1}-\tilde{y}_i^{k+1} \rangle-(1-\sigma)\langle x^\star-\tilde{y}_i^{k+1}, \nabla g_i(x^{k+1}) \rangle,~i\in[m-1].
\end{aligned}
\end{equation*}
By summing the second inequality from $i=1$ to $m-1$ and adding it to the first inequality, we obtain
\begin{equation}\label{t1}
\begin{aligned}
F(x^\star)\ge& \sum_{i=1}^{m-1}f_i(\tilde{y}_i^{k+1})+f_m(x^{k+1})+\frac{1}{\gamma}\sum_{i=1}^{m-1}\langle x^\star-z_i^k+x^{k+1}-\tilde{y}_i^{k+1},x^{k+1}-\tilde{y}_i^{k+1} \rangle\\
&+\frac{1}{\gamma}\sum_{i=1}^{m-1}\alpha_i\left(\langle x^\star-\tilde{y}_i^{k+1}, x^{k+1}-\tilde{y}_i^{k+1} \rangle+\langle x^\star-x^{k+1}, y_i^k-x^{k+1} \rangle\right)\\
&-\frac{1}{m-1}\sum_{i=1}^{m-1}\langle x^\star-x^{k+1}, \nabla g_m(y_i^k) \rangle-\sigma\sum_{i=1}^{m-1}\langle x^\star-x^{k+1}, \nabla g_i(y_i^k) \rangle\\
&-(1-\sigma)\sum_{i=1}^{m-1}\langle x^\star-\tilde{y}_i^{k+1}, \nabla g_i(x^{k+1}) \rangle\\
=& \sum_{i=1}^{m-1}f_i(\tilde{y}_i^{k+1})+f_m(x^{k+1})+\frac{1}{\gamma}\sum_{i=1}^{m-1}\langle x^\star-z_i^k, x^{k+1}-\tilde{y}_i^{k+1} \rangle+\frac{1}{\gamma}\sum_{i=1}^{m-1}\|x^{k+1}-\tilde{y}_i^{k+1}\|^2\\
&+\frac{1}{2\gamma}\sum_{i=1}^{m-1}\alpha_i\left(\|\tilde{y}_i^{k+1}-x^\star\|^2-\|y^k_i-x^\star\|^2+\|x^{k+1}-\tilde{y}_i^{k+1}\|^2+\|x^{k+1}-y_i^k\|^2 \right)\\
&-\frac{1}{m-1}\sum_{i=1}^{m-1}\langle x^\star-x^{k+1}, \nabla g_m(y_i^k) \rangle-\sigma\sum_{i=1}^{m-1}\langle x^\star-x^{k+1}, \nabla g_i(y_i^k) \rangle\\
&-(1-\sigma)\sum_{i=1}^{m-1}\langle x^\star-\tilde{y}_i^{k+1}, \nabla g_i(x^{k+1}) \rangle.
\end{aligned}
\end{equation}

Next, we estimate the remaining cross terms in \eqref{t1}. Using \eqref{all4}, we obtain
\begin{equation}\label{t2}
\begin{aligned}
&\langle x^\star-z_i^k, x^{k+1}-\tilde{y}_i^{k+1} \rangle\\
=&\langle x^\star-z_i^\star, x^{k+1}-\tilde{y}_i^{k+1} \rangle+\langle z_i^\star-z_i^k, x^{k+1}-\tilde{y}_i^{k+1} \rangle\\
=&\langle x^\star-z_i^\star, x^{k+1}-\tilde{y}_i^{k+1} \rangle+\frac{1}{\lambda_i}\langle z_i^\star-z_i^k, z_i^{k}-\tilde{z}_i^{k+1} \rangle\\
=&\langle x^\star-z_i^\star, x^{k+1}-\tilde{y}_i^{k+1} \rangle+\frac{1}{2\lambda_i}\left( \|\tilde{z}_i^{k+1}-z_i^\star\|^2-\|{z}_i^{k}-z_i^\star\|^2-\|\tilde{z}_i^{k+1}-z_i^k\|^2\right),
\end{aligned}
\end{equation}
for all $i\in [m-1]$. Letting $g=g_m$, $y=x^{k+1}$, $z=x^\star$, $x=y_i^k$, and $L=1/\beta_m$ in $\eqref{dc}$, we derive
\begin{equation}\label{d1}
  \langle x^\star-x^{k+1}, \nabla g_m(y_i^k) \rangle \le g_m(x^\star) - g_m(x^{k+1}) + \frac{1}{2\beta_m}\|x^{k+1} - y_i^k\|^2, \quad i \in [m-1].
\end{equation}
Similarly, we derive
\begin{equation}\label{d2}
  \langle x^\star-x^{k+1}, \nabla g_i(y_i^k) \rangle \le g_i(x^\star) - g_i(x^{k+1}) + \frac{1}{2\beta_i}\|x^{k+1} - y_i^k\|^2, \quad i \in [m-1],
\end{equation}
and
\begin{equation}\label{d3}
  \langle x^\star - \tilde{y}_i^{k+1}, \nabla g_i(x^{k+1}) \rangle \le g_i(x^\star) - g_i(\tilde{y}_i^{k+1}) + \frac{1}{2\beta_i}\|x^{k+1} - \tilde{y}_i^{k+1}\|^2, \quad i \in [m-1].
\end{equation}
Substituting \eqref{t2}--\eqref{d3} into \eqref{t1} and simplifying yields
\begin{equation}\label{t3}
\begin{aligned}
F(x^\star)\ge& \sum_{i=1}^{m-1}\left(f_i(\tilde{y}_i^{k+1})+\sigma g_i(x^{k+1})+(1-\sigma) g_i(\tilde{y}_i^{k+1})\right)+f_m(x^{k+1})+g_m(x^{k+1})\\
&+\frac{1}{2\gamma}\sum_{i=1}^{m-1}\left(\alpha_i\|\tilde{y}_i^{k+1}-x^\star\|^2+\frac{1}{\lambda_i}\|\tilde{z}_i^{k+1}-z_i^\star\|^2-\alpha_i\|y^k_i-x^\star\|^2-\frac{1}{\lambda_i}\|{z}_i^{k}-z_i^\star\|^2\right)\\
&+\frac{1}{2\gamma}\sum_{i=1}^{m-1}\left(2+\alpha_i-\lambda_i-\frac{(1-\sigma)\gamma}{\beta_i}\right)\|x^{k+1}-\tilde{y}_i^{k+1}\|^2\\
&+\frac{1}{\gamma}\sum_{i=1}^{m-1}\langle x^\star-z_i^\star, x^{k+1}-\tilde{y}_i^{k+1} \rangle\\
&+\frac{1}{2\gamma}\sum_{i=1}^{m-1}\left(\alpha_i-\frac{\gamma}{(m-1)\beta_m}-\frac{\sigma\gamma}{\beta_i}\right)\|x^{k+1}-y_i^k\|^2.
\end{aligned}
\end{equation}
Further substituting \eqref{e3} and \eqref{e4} into \eqref{t3} implies
\begin{equation}\label{t4}
\begin{aligned}
F(x^\star)\ge& \sum_{i=1}^{m-1}\left(f_i(\tilde{y}_i^{k+1})+\sigma g_i(x^{k+1})+(1-\sigma) g_i(\tilde{y}_i^{k+1})\right)+f_m(x^{k+1})+g_m(x^{k+1})\\
&+\frac{1}{2\gamma}\sum_{i=1}^{m-1}\left\{\frac{\alpha_i}{p_i}\mathbb{E}\left(\|{y}_i^{k+1}-x^\star\|^2\mid\mathcal{F}_k\right)+\frac{1}{\lambda_ip_i}\mathbb{E}\left(\|{z}_i^{k+1}-z_i^\star\|^2\mid\mathcal{F}_k\right)\right\}\\
&-\frac{1}{2\gamma}\sum_{i=1}^{m-1}\left(\frac{\alpha_i}{p_i}\|y^k_i-x^\star\|^2+\frac{1}{\lambda_ip_i}\|{z}_i^{k}-z_i^\star\|^2\right)+\frac{1}{\gamma}\sum_{i=1}^{m-1}\langle x^\star-z_i^\star, x^{k+1}-\tilde{y}_i^{k+1} \rangle\\
&+\frac{1}{2\gamma}\sum_{i=1}^{m-1}\left(2+\alpha_i-\lambda_i-\frac{(1-\sigma)\gamma}{\beta_i}\right)\|x^{k+1}-\tilde{y}_i^{k+1}\|^2\\
&+\frac{1}{2\gamma}\sum_{i=1}^{m-1}\left(\alpha_i-\frac{\gamma}{(m-1)\beta_m}-\frac{\sigma\gamma}{\beta_i}\right)\|x^{k+1}-y_i^k\|^2.
\end{aligned}
\end{equation}
The proof is completed.

}
\paragraph{ The proof of Theorem 2 (Rate).}{\rm (i): Since $\mathbb{E}\|z_i^{k+1}\|<\infty$ for each $i \in [m-1]$, the conditional expectation $\mathbb{E}(z^{k+1}_i \mid \mathcal{F}_k)$ exists, and 
\[
\mathbb{E}(z^{k+1}_i \mid \mathcal{F}_k) = \sum_{j=1}^n \mathbb{E}\left( \langle z^{k+1}_i, e^j \rangle \big| \mathcal{F}_k \right) e^j, \quad i \in [m-1],
\]
where $\{e^j\}_{j=1}^n$ denotes the standard orthonormal basis of $\mathbb{R}^n$. By an argument analogous to that in the proof of Lemma \ref{ce}, it follows that
\[
\mathbb{E}\left( \langle z^{k+1}_i, e^j \rangle \big| \mathcal{F}_k \right) = p_i \langle \tilde{z}^{k+1}_i, e^j \rangle + (1-p_i) \langle z^{k}_i, e^j \rangle, \quad \forall j \in [n].
\]
Applying Parseval's equality yields
\begin{equation}\label{zce}
    \mathbb{E}\left( z^{k+1}_i - z_i^k \big| \mathcal{F}_k \right) = p_i \left( \tilde{z}^{k+1}_i - z_i^k \right) = p_i \lambda_i \left( \tilde{y}_i^{k+1} - x^{k+1} \right), \quad i \in [m-1].
\end{equation}
For any $K \in \mathbb{N}$ with $K \geq 1$, define the averages
\[
x_{\mathrm{av}}^K := \frac{1}{K} \sum_{k=0}^{K-1} x^{k+1}, \quad \tilde{y}_{\mathrm{av},i}^K := \frac{1}{K} \sum_{k=0}^{K-1} \tilde{y}_i^{k+1}, \quad i \in [m-1].
\]
Then, for each $i \in [m-1]$, it holds that
\begin{equation}\label{ere1}
\begin{aligned}
    \left\| \mathbb{E}\left( x_{\mathrm{av}}^K - \tilde{y}_{\mathrm{av},i}^K \right) \right\| 
    &= \frac{1}{K} \left\| \mathbb{E} \left[ \sum_{k=0}^{K-1} \left( x^{k+1} - \tilde{y}_i^{k+1} \right) \right] \right\| \\
    &= \frac{1}{\lambda_i p_i K} \left\| \sum_{k=0}^{K-1} \mathbb{E} \left[ \mathbb{E} \left( z_i^{k+1} - z_i^k \mid \mathcal{F}_k \right) \right] \right\| \\
    &= \frac{1}{\lambda_i p_i K} \left\| \sum_{k=0}^{K-1} \mathbb{E} \left( z_i^{k+1} - z_i^k \right) \right\| \\
    &\leq \frac{1}{\lambda_i p_i K} \mathbb{E} \left( \| z_i^K - z_i^0 \| \right).
\end{aligned}
\end{equation}
Since the sequence $\{z_i^k\}_{k=0}^\infty$ is bounded for each $i \in [m-1]$, it follows that
\begin{equation}\label{cr}
\left\| \mathbb{E}\left( x_{\mathrm{av}}^K - \tilde{y}_{\mathrm{av},i}^K \right) \right\| \leq \frac{\sup_{i,K} \{ \| z_i^K - z_i^0 \| \} }{ \inf_i \{ \lambda_i p_i \} K }, \quad \forall i \in [m-1].
\end{equation}
By employing the analogous proof technique, we similarly obtain
\begin{equation}\label{cr2}
\left\| \mathbb{E}\left( y_{\mathrm{av},i}^K - \tilde{y}_{\mathrm{av},i}^K \right) \right\| \leq \frac{\sup_{i,K} \{ \| y_i^K - y_i^0 \| \} }{ \inf_i \{ p_i \} K }, \quad \forall i \in [m-1].
\end{equation}

(ii). Summing both sides of \eqref{t4} over $k=0$ to $K-1$ yields
\begin{equation}\label{ef}
\begin{aligned}
& \sum_{k=0}^{K-1} \mathbb{E} \left\{ \sum_{i=1}^{m-1} \left(f_i(\tilde{y}_i^{k+1}) + \sigma g_i(x^{k+1}) + (1-\sigma) g_i(\tilde{y}_i^{k+1}) \right) + f_m(x^{k+1}) + g_m(x^{k+1}) \right\} \\
\leq  & K F(x^\star) + \frac{1}{2\gamma} (a_0 - a_K) - \frac{1}{2\gamma} \sum_{k=0}^{K-1} \sum_{i=1}^{m-1} \left( 2 + \alpha_i - \lambda_i - \frac{(1-\sigma) \gamma}{\beta_i} \right) \mathbb{E} \left( \| \tilde{y}_i^{k+1} - x^{k+1} \|^2 \right) \\
&- \frac{1}{2\gamma} \sum_{k=0}^{K-1} \sum_{i=1}^{m-1} \left( \alpha_i - \frac{\gamma}{(m-1) \beta_m} - \frac{\sigma \gamma}{\beta_i} \right) \mathbb{E} \left( \| x^{k+1} - y_i^k \|^2 \right) \\
& + \frac{1}{\gamma} \sum_{i=1}^{m-1} \mathbb{E} \left[ \left\langle x^\star - z_i^\star, \sum_{k=0}^{K-1} \left( x^{k+1} - \tilde{y}_i^{k+1} \right) \right\rangle \right],
\end{aligned}
\end{equation}

Applying Jensen's inequality yields
\begin{equation}\label{ef1}
\begin{aligned}
& \mathbb{E} \left\{ \sum_{i=1}^{m-1} \left( f_i(\tilde{y}_{\mathrm{av},i}^K) + \sigma g_i(x_{\mathrm{av}}^K) + (1-\sigma) g_i(\tilde{y}_{\mathrm{av},i}^K) \right) + f_m(x_{\mathrm{av}}^K) + g_m(x_{\mathrm{av}}^K) \right\} - F(x^\star) \\
\leq {} & \frac{a_0}{2 \gamma K} - \frac{1}{2 \gamma K} \sum_{k=0}^{K-1} \sum_{i=1}^{m-1} \left( 2 + \alpha_i - \lambda_i - \frac{(1-\sigma) \gamma}{\beta_i} \right) \mathbb{E} \left( \| \tilde{y}_i^{k+1} - x^{k+1} \|^2 \right) \\
& - \frac{1}{2 \gamma K} \sum_{k=0}^{K-1} \sum_{i=1}^{m-1} \left( \alpha_i - \frac{\gamma}{(m-1) \beta_m} - \frac{\sigma \gamma}{\beta_i} \right) \mathbb{E} \left( \| x^{k+1} - y_i^k \|^2 \right) \\
& + \frac{1}{\gamma} \sum_{i=1}^{m-1} \mathbb{E} \left[ \left\langle x^\star - z_i^\star, x_{\mathrm{av}}^K - \tilde{y}_{\mathrm{av},i}^K \right\rangle \right].
\end{aligned}
\end{equation}

By the definition of expectation and the Cauchy-Schwarz inequality, it holds that
\begin{equation}\label{ef2}
\begin{aligned}
\Big| \mathbb{E} \left[ \left\langle x^\star - z_i^\star, x_{\mathrm{av}}^K - \tilde{y}_{\mathrm{av},i}^K \right\rangle \right] \Big| &= \Big| \left\langle x^\star - z_i^\star, \mathbb{E} \left( x_{\mathrm{av}}^K - \tilde{y}_{\mathrm{av},i}^K \right) \right\rangle \Big| \\
&\leq \| x^\star - z_i^\star \| \left\| \mathbb{E} \left( x_{\mathrm{av}}^K - \tilde{y}_{\mathrm{av},i}^K \right) \right\|, \quad \forall i \in [m-1].
\end{aligned}
\end{equation}

Substituting \eqref{cr} and \eqref{ef2} into \eqref{ef1} yields
\begin{equation}\label{ef7}
\begin{aligned}
& \mathbb{E} \left\{ \sum_{i=1}^{m-1} \left(f_i(\tilde{y}_{\mathrm{av},i}^K) +\sigma g_i(x_{\mathrm{av}}^K) + (1-\sigma) g_i(\tilde{y}_{\mathrm{av},i}^K) \right) + f_m(x_{\mathrm{av}}^K) + g_m(x_{\mathrm{av}}^K) \right\} - F(x^\star) \\
\leq {} & \frac{a_0}{2 \gamma K} - \frac{1}{2 \gamma K} \sum_{k=0}^{K-1} \sum_{i=1}^{m-1} \left[ 2 + \alpha_i - \lambda_i - \frac{(1-\sigma) \gamma}{\beta_i} \right] \mathbb{E} \left( \| \tilde{y}_i^{k+1} - x^{k+1} \|^2 \right) \\
& - \frac{1}{2 \gamma K} \sum_{k=0}^{K-1} \sum_{i=1}^{m-1} \left( \alpha_i - \frac{\gamma}{(m-1) \beta_m} - \frac{\sigma \gamma}{\beta_i} \right) \mathbb{E} \left( \| x^{k+1} - y_i^k \|^2 \right) + \frac{c_1}{K},
\end{aligned}
\end{equation}
where $$c_1 = \frac{(m-1) \sup_{K,i} \{ \| z_i^\star - x^\star \| \cdot \| z_i^K - z_i^0 \| \}}{\gamma \inf_i \{ \lambda_i p_i \}}.$$
Define $c_2 = \min_{i \in [m-1]} \left\{ 2 + \alpha_i - \frac{(1-\sigma)\gamma}{2\beta_i} - \lambda_i, \quad \alpha_i - \frac{\gamma}{2\beta_m (m-1)} - \frac{\sigma \gamma}{2 \beta_i} \right\}$, which satisfies $c_2 > 0$. Rearranging the terms in \eqref{m41} yields the inequality
\begin{equation}\label{e7}
\mathbb{E}\left[ \sum_{i=1}^{m-1} \left( \| x^{k+1} - \tilde{y}_i^{k+1} \|^2 + \| x^{k+1} - y_i^k \|^2 \right) \right] \leq \frac{1}{c_2} \left( a_k - a_{k+1} \right),
\end{equation}
where   
$$
a_k=\sum_{i=1}^{m-1}\mathbb{E}\left(\frac{\alpha_i}{p_i}\|y_i^k-x^*\|^2+\frac{1}{\lambda p_i}\|z_i^k-z_i^*\|^2\right).
$$
Summing \eqref{e7} over $k=0$ to $K-1$ gives
\begin{equation}\label{rr2}
\sum_{k=0}^{K-1} \sum_{i=1}^{m-1} \mathbb{E} \left( \| x^{k+1} - \tilde{y}_i^{k+1} \|^2 + \| x^{k+1} - y_i^k \|^2 \right) \leq \frac{a_0}{c_2}.
\end{equation}
If $\gamma \in \left( 0, \min_{i \in [m-1]} \left\{ \frac{\alpha_i}{\frac{1}{\beta_m (m-1)} + \frac{\sigma}{\beta_i}}, \frac{(2 + \alpha_i) \beta_i}{1-\sigma} \right\} \right]$
and $\lambda_i \in \left( 0, 2 + \alpha_i - \frac{(1-\sigma) \gamma}{\beta_i} \right]$, then the following inequalities hold:
\[
\begin{cases}
\alpha_i - \frac{\gamma}{(m-1) \beta_m} - \frac{\sigma \gamma}{\beta_i} \geq 0, \\
2 + \alpha_i - \lambda_i - \frac{(1-\sigma) \gamma}{\beta_i} \geq 0.
\end{cases}
\]
In this case, \eqref{ef7} implies that
\begin{equation}\label{ef8}
\begin{aligned}
& \mathbb{E} \left\{ \sum_{i=1}^{m-1} \left( f_i(\tilde{y}_{\mathrm{av},i}^K) + \sigma g_i(x_{\mathrm{av}}^K) + (1-\sigma) g_i(\tilde{y}_{\mathrm{av},i}^K) \right) + f_m(x_{\mathrm{av}}^K) + g_m(x_{\mathrm{av}}^K) \right\} - F(x^\star) \\
\leq {} & \frac{a_0}{2 \gamma K} + \frac{c_1}{K}.
\end{aligned}
\end{equation}
Otherwise, the following inequalities hold:
\[
\begin{cases}
\alpha_i - \frac{\gamma}{(m-1) \beta_m} - \frac{\sigma \gamma}{\beta_i} < 0, \\
2 + \alpha_i - \lambda_i - \frac{(1-\sigma) \gamma}{\beta_i} < 0.
\end{cases}
\]
Let \( c_3 \) be defined by
\[
c_3 = \max_{i \in [m-1]} \left\{ \frac{\gamma}{(m-1) \beta_m} + \frac{\sigma \gamma}{\beta_i} - \alpha_i, \quad \lambda_i + \frac{(1-\sigma) \gamma}{\beta_i} - 2 - \alpha_i \right\},
\]
which satisfies \( c_3 > 0 \).
Under these conditions, combining \eqref{rr2} and \eqref{ef7} yields
\begin{equation}\label{ef9}
\begin{aligned}
& \mathbb{E} \left\{ \sum_{i=1}^{m-1} \left( f_i(\tilde{y}_{\mathrm{av},i}^K) + \sigma g_i(x_{\mathrm{av}}^K) + (1-\sigma) g_i(\tilde{y}_{\mathrm{av},i}^K) \right) + f_m(x_{\mathrm{av}}^K) + g_m(x_{\mathrm{av}}^K) \right\} - F(x^\star) \\
\leq{} & \frac{a_0}{2 \gamma K} + \frac{c_3 a_0}{2 c_2 \gamma K} + \frac{c_1}{K}.
\end{aligned}
\end{equation}
For each \( i \in [m-1] \), set \( x_i = \tilde{y}_{\mathrm{av},i}^K \) and \( x_m = x_{\mathrm{av}}^K \) in \eqref{fc1}. Then, applying the Cauchy-Schwarz inequality yields
\begin{equation}\label{ef10}
\begin{aligned}
& \mathbb{E} \left\{ \sum_{i=1}^m f_i(\tilde{y}_{\mathrm{av},i}^K) + \sum_{i=1}^{m-1} \left(\sigma g_i(x_{\mathrm{av}}^K)+ (1-\sigma) g_i(\tilde{y}_{\mathrm{av},i}^K)  \right) + g_m(x_{\mathrm{av}}^K) \right\} - F(x^\star) \\
\geq{} & \frac{1}{\gamma} \sum_{i=1}^{m-1} \mathbb{E} \left[ \langle x^\star - z_i^\star, x_{\mathrm{av}}^K - \tilde{y}_{\mathrm{av},i}^K \rangle \right] \\
\geq{} & - \frac{1}{\gamma} \sum_{i=1}^{m-1} \| x^\star - z_i^\star \| \left\| \mathbb{E} \left( x_{\mathrm{av}}^K - \tilde{y}_{\mathrm{av},i}^K \right) \right\| \\
\geq{} & - \frac{c_1}{K} \sum_{i=1}^{m-1} \frac{\| x^\star - z_i^\star \|}{\gamma}.
\end{aligned}
\end{equation}
Therefore, by combining \eqref{ef8}, \eqref{ef9}, and \eqref{ef10}, we deduce that
\[
\mathbb{E} \left\{ \sum_{i=1}^m f_i(\tilde{y}_{\mathrm{av},i}^K) + \sum_{i=1}^{m-1} \left(\sigma g_i(x_{\mathrm{av}}^K)+ (1-\sigma) g_i(\tilde{y}_{\mathrm{av},i}^K)  \right) + g_m(x_{\mathrm{av}}^K) \right\} - F(x^\star) = O\left( \frac{1}{K} \right).
\]
For each \( i \in [m-1] \), it can be shown that
\begin{equation}\label{yk}
\begin{cases}
\mathbb{E} f_i(\tilde{y}_i^{k+1}) = \frac{1}{p_i} \mathbb{E} f_i(y_i^{k+1}) - \frac{1}{p_i} \mathbb{E} f_i(y_i^k) + \mathbb{E} f_i(y_i^k), \\
\mathbb{E} g_i(\tilde{y}_i^{k+1}) = \frac{1}{p_i} \mathbb{E} g_i(y_i^{k+1}) - \frac{1}{p_i} \mathbb{E} g_i(y_i^k) + \mathbb{E} g_i(y_i^k).
\end{cases}
\end{equation}
By substituting \eqref{yk} into Lemma 4 and repeating the preceding arguments, we obtain
\begin{equation}
\begin{aligned}
& \mathbb{E} \left\{ \sum_{i=1}^{m-1} \left( f_i(y_{\mathrm{av},i}^K) +\sigma g_i(x_{\mathrm{av}}^K) + (1-\sigma) g_i(y_{\mathrm{av},i}^K) \right) + f_m(x_{\mathrm{av}}^K) + g_m(x_{\mathrm{av}}^K) \right\} - F(x^\star) \\
\leq{} & \frac{a_0}{2 \gamma K} + \frac{c_3 a_0}{2 c_2 \gamma K} + \frac{c_1}{K} + \frac{ \sum_{i=1}^{m-1} \frac{1}{p_i} \left[ \mathbb{E} f_i(y_i^0) + (1-\sigma) \mathbb{E} g_i(y_i^0) - \mathbb{E} f_i(y_i^K) - (1-\sigma) \mathbb{E} g_i(y_i^K) \right] }{K}.
\end{aligned}
\end{equation}
Since the sequences \(\{ y_i^k \}_{k=1}^\infty\) are bounded for each \( i \in [m-1] \) and the functions \( f_i, g_i \) are lower semicontinuous, there exists a positive constant \( c_4 \) satisfying
\[
\sum_{i=1}^{m-1} \frac{1}{p_i} \Big| \mathbb{E} f_i(y_i^0) + (1-\sigma) \mathbb{E} g_i(y_i^0) - \mathbb{E} f_i(y_i^K) - (1-\sigma) \mathbb{E} g_i(y_i^K) \Big| \leq c_4, \quad \forall K \geq 1.
\]
For each \( i \in [m-1] \), set \( x_i = y_{\mathrm{av},i}^K \) and \( x_m = x_{\mathrm{av}}^K \) in \eqref{fc1}. Then, by applying the Cauchy-Schwarz inequality, we have
\begin{equation}\label{ef11}
\begin{aligned}
&  \mathbb{E} \left\{ \sum_{i=1}^{m-1} \left( f_i(y_{\mathrm{av},i}^K) +\sigma g_i(x_{\mathrm{av}}^K) + (1-\sigma) g_i(y_{\mathrm{av},i}^K) \right) + f_m(x_{\mathrm{av}}^K) + g_m(x_{\mathrm{av}}^K) \right\} - F(x^\star) \\
\geq{} & \frac{1}{\gamma} \sum_{i=1}^{m-1} \mathbb{E} \left[ \langle x^\star - z_i^\star, x_{\mathrm{av}}^K - y_{\mathrm{av},i}^K \rangle \right] \\
\geq{} & - \frac{1}{\gamma} \sum_{i=1}^{m-1} \| x^\star - z_i^\star \| \left\| \mathbb{E} \left( x_{\mathrm{av}}^K - y_{\mathrm{av},i}^K \right) \right\| \\
\geq{} & - \frac{1}{\gamma} \sum_{i=1}^{m-1} \| x^\star - z_i^\star \| \left( \left\| \mathbb{E} \left( x_{\mathrm{av}}^K - \tilde{y}_{\mathrm{av},i}^K \right) \right\| + \left\| \mathbb{E} \left( y_{\mathrm{av},i}^K - \tilde{y}_{\mathrm{av},i}^K \right) \right\| \right) \\
\geq{} & - \frac{c_5}{K} \sum_{i=1}^{m-1} \frac{ \| x^\star - z_i^\star \| }{\gamma},
\end{aligned}
\end{equation}
where 
\[
c_5 = \frac{ \sup_{i,K} \| z_i^K - z_i^0 \| }{ \inf_i \{ \lambda_i p_i \} } + \frac{ \sup_{i,K} \| y_i^K - y_i^0 \| }{ \inf_i \{ p_i \} }.
\]
This completes the proof.

}

\end{document}